\documentclass[paper=a4, USenglish, numbers=noenddot]{scrartcl} 
\pdfoutput=1

\usepackage{cmap}		
\usepackage[utf8]{inputenc}	
\usepackage[T1]{fontenc}	
\usepackage[USenglish]{babel}

\usepackage{microtype}

\usepackage{amsmath}
\usepackage{amssymb}
\usepackage{mathtools}		

\usepackage{tikz}
\usetikzlibrary{intersections,calc,through}

\usepackage[pdftitle={Wild translation surfaces and infinite genus},pdfauthor={Anja Randecker},pdfborder={0 0 0}]{hyperref}

\usepackage[figure]{hypcap}
\usepackage{enumitem}

\usepackage{aliascnt}
\usepackage{amsthm}

\newtheoremstyle{break}
 {} 
 {} 
 {\itshape} 
 {} 
 {\bfseries} 
 {} 
 {\newline} 
 {\thmname{#1}\thmnumber{ #2}\thmnote{ (#3)}} 
 
\newtheoremstyle{breakdef}
 {} 
 {} 
 {} 
 {} 
 {\bfseries} 
 {} 
 {\newline} 
 {\thmname{#1}\thmnumber{ #2}\thmnote{ (#3)}} 
 
\newtheoremstyle{remark}
 {} 
 {} 
 {} 
 {} 
 {\itshape} 
 {.} 
 {0.5em} 
 {\thmname{#1}\thmnumber{ #2}\thmnote{ {\normalfont (}#3{\normalfont )}}} 

\theoremstyle{breakdef}
\newtheorem{definition}{Definition}[section]

\theoremstyle{remark}

\newaliascnt{rem}{definition}  
\newtheorem{rem}[rem]{Remark}
\aliascntresetthe{rem}

\newaliascnt{exa}{definition}  
\newtheorem{exa}[exa]{Example}
\aliascntresetthe{exa} 

\newaliascnt{lem}{definition}  
\newtheorem{lem}[lem]{Lemma}
\aliascntresetthe{lem}

\newaliascnt{conv}{definition}  
\newtheorem{conv}[conv]{Convention}
\aliascntresetthe{conv}

\theoremstyle{break}

\newtheorem{thm}{Theorem}

\newaliascnt{prop}{definition}  
\newtheorem{prop}[prop]{Proposition}
\aliascntresetthe{prop}

\newaliascnt{cor}{definition}  
\newtheorem{cor}[cor]{Corollary}
\aliascntresetthe{cor}

\DeclareMathOperator{\ir}{ir}
\DeclareMathOperator{\im}{im}
\DeclareMathOperator{\ord}{ord}

\renewcommand{\Re}{\mathrm{Re}}
\renewcommand{\Im}{\mathrm{Im}}

\author{Anja Randecker}
\title{Wild translation surfaces \\ and infinite genus}
\date{\normalsize{
University of Toronto\\
Toronto, Canada\\
e-mail: \href{mailto:anja@math.toronto.edu}{anja@math.toronto.edu}
}}

\begin{document}

\def\sectionautorefname{Section}
 
\maketitle

\begin{abstract}
\noindent
The Gauss--Bonnet formula for classical translation surfaces relates the cone angle of the singularities (geometry) to the genus of the surface (topology).
When considering more general translation surfaces, we observe so-called wild singularities for which the notion of cone angle is not applicable any more.

We study whether there still exist relations between the geometry and the topology for translation surfaces with wild singularities. By considering short saddle connections, we determine under which conditions the existence of a wild singularity implies infinite genus. We apply this to show that parabolic or essentially finite translation surfaces with wild singularities have infinite~genus.
\end{abstract}

Classical translation surfaces are objects at the intersection of many different fields, such as dynamical systems, Teichmüller theory, algebraic geometry, topology, and geometric group theory.
The history of translation surfaces starts in the time of the article \cite{fox_kershner_36}. Fox and Kershner obtained translation surfaces in the theory of billiards when ``unfolding'' polygons with rational angles.

The most visual description of classical translation surfaces is given by considering finitely many polygons in the Euclidean plane.
If every edge of the polygons can be identified with a parallel edge of the same length so that we obtain a connected, orientable surface then the resulting object is a translation surface. It is locally flat at all points with the possible exception of the former vertices of the polygons. These exceptional points are called \emph{singularities} and they are cone points of the resulting surface with cone angle $2\pi k$ for some $k\geq 2$.

A natural generalization is to drop the condition that the number of polygons has to be finite. When gluing infinitely many polygons, the local flatness still holds but the behaviour of the singularities is more diverse than in the classical case. This kind of translation surface is often called \emph{infinite} in the literature, but we will not follow this convention here and simply call it a \emph{translation surface}.

Recently, the interest in this generalization of translation surfaces has grown:
There are results on Veech groups by Chamanara \cite{chamanara_04}, Hubert and Schmithüsen \cite{hubert_schmithuesen_10}, and Przytycki, Schmithüsen, and Valdez \cite{przytycki_schmithuesen_valdez_11}, results on the dynamics by Hooper \cite{hooper_14a}, Treviño \cite{trevino_12}, and Lindsey and Treviño \cite{lindsey_trevino_16}, and results on infinite coverings of finite translation surfaces (especially for the wind-tree model) by Delecroix, Hubert, and Lelièvre \cite{delecroix_hubert_lelievre_11}, Hubert, Lelièvre, and Troubetzkoy \cite{hubert_lelievre_troubetzkoy_11}, Avila and Hubert \cite{avila_hubert_12}, Hooper and Weiss \cite{hooper_weiss_12}, Hooper, Hubert, and Weiss \cite{hooper_hubert_weiss_13}, and Fr{\c{a}}czek and Ulcigrai \cite{fraczek_ulcigrai_14}.
However, while we have a classification of finite translation surfaces by studying strata of the moduli space, there is no systematic description for the generalized ones so far. A first step towards such a classification can be to understand and classify the singularities of the translation surfaces.

When considering translation surfaces with interesting singularities, it is natural to start with translation surfaces with exactly one singularity. So, a lot of the recently described examples are Loch Ness monsters, i.e.\ surfaces with infinite genus and one end (see Richards \cite{richards_63} for the definition of ends and the classification of noncompact surfaces using ends). For instance, these examples include the baker's map surface by Chamanara \cite{chamanara_04}, the Arnoux--Yoccoz surfaces as in Bowman \cite{bowman_13}, the ladder surface by Hooper \cite{hooper_15}, the stack of boxes by Bowman \cite{bowman_12}, and the exponential surface in the author's thesis \cite{randecker_16}.

For all these mentioned translation surfaces, the singularity is a so-called wild singularity and the topology is always the same, namely infinite genus and one end. The number of ends can easily be increased by gluing in half-cylinders but it is hard to reduce the genus. This leads to the conjecture that the existence of a wild singularity implies infinite~genus.

In full generality, this conjecture is not true, as we will see in \autoref{exa_nested_cylinders} but in our main theorem we provide necessary conditions for a wild singularity to imply infinite~genus.
These conditions include \emph{xossiness} (short for e\emph{x}istence \emph{o}f \emph{s}hort \emph{s}addle connections \emph{i}ntersected \emph{n}ot by \emph{e}ven \emph{s}horter \emph{s}addle connections) which describes the existence of arbitrarily short geodesic segments between singularities so that there is a lower bound on the length of intersecting segments of that type. We motivate in \autoref{sec_xossiness} that xossiness is not a strong condition, as it is fulfilled, for instance, if the geodesic flow is defined for all times in almost every direction.

\begin{thm}
 Let $(X,\mathcal{A})$ be a translation surface with the following properties:
 \begin{enumerate}
  \item The singularities of $(X,\mathcal{A})$ are discrete.
  \item There exists a wild singularity $\sigma$ that fulfills xossiness.
  \item There exist two directions $\theta_1, \theta_2$ for which the geodesic flows $F_{\theta_1}, F_{\theta_2}$ are recurrent.
 \end{enumerate}
 Then $X$ has infinite genus.
\end{thm}
\setcounter{thm}{0}

\noindent
For example, this theorem is applicable for all of the Loch Ness monster examples listed above (although for each of them it is easy to check by hand that the genus is infinite).

On the other hand, the theorem is not applicable to so-called \emph{half-translation surfaces} or \emph{reflection-translation surfaces}, where edges of the polygons can be identified by translations or by reflections. In particular, there exist infinite reflection-translation surfaces with wild singularities and genus $0$ such as the ones coming from the construction by de Carvalho and Hall \cite[Theorem 42]{carvalho_hall_12}. This might seem contradictory as it is a standard procedure to pass from a reflection-translation surface to a double cover which is a translation surface to use the richer theory of translation surfaces.
However, when doubling a capped-off cylinder in a reflection-translation surface, it might form a handle. This introduces additional genus and so the double cover of a genus $0$ reflection-translation surface can indeed be of infinite genus.

The article is structured as follows. We start with the necessary definitions on translation surfaces, their singularities, and rotational components in \autoref{sec_definitions}.
In \autoref{sec_criterion_infinite_genus}, a criterion is developed how infinite genus can be detected using saddle connections, i.e.\ geodesics that connect singularities.
We establish the existence of short saddle connections in \autoref{sec_saddle_connections_intersections} and study the existence of lower bounds of the length of intersecting saddle connections.
In \autoref{sec_xossiness}, we define the notion of xossiness and illustrate it by proving two criteria that imply xossiness.
The proof of the main theorem is carried out in \autoref{sec_theorem_corollaries_generalizations}, followed by two applications and a discussion on whether the conditions in the main theorem are necessary.

\bigskip
\noindent\textbf{Acknowledgments.} This article was written as part of the author's PhD project. The author wants to thank her advisor Gabi Weitze-Schmithüsen for hundreds of valuable discussions and her constant encouragement, Pat Hooper for many helpful discussions (especially at the beginning of this project, see \autoref{exa_nested_cylinders}), and Ferrán Valdez for sharing his thoughts on singularities of translation surfaces, starting with explanations on \cite{bowman_valdez_13}. Furthermore, the author wants to thank Joshua Bowman, Vincent Delecroix, Rodrigo Treviño, Chenxi Wu, Lucien Clavier, and Rodolfo Rios-Zertuche for helpful discussions and hints.

The author was partially supported within the project ``Dynamik unendlicher Translationsflächen'' in the ``Juniorprofessoren-Programm'' of the Ministry of Science, Research and the Arts.

\section{Basics on translation surfaces and their singularities} \label{sec_definitions}

As there are several nuances in generalizing the notion of classical translation surfaces, we first make clear the definition that we use in this article.

\begin{definition}[(Finite) translation surface] \label{def_finite_translation_surfaces}
 A \emph{translation surface} $(X,\mathcal{A})$ is a connected surface $X$ with a translation structure $\mathcal{A}$ on~$X$, i.e.\ a maximal atlas on $X$ such that the transition functions are translations.
 Via the translation structure we can pull back the Euclidean metric from $\mathbb{R}^2$ to $X$.
 
 The translation surface $(X,\mathcal{A})$ is called \emph{finite} if the metric completion $\overline{X}$ is a compact surface and $\overline{X}\setminus X$ is discrete.
\end{definition}

\noindent
Note that there exist examples of translation surfaces which are not finite and such that the metric completion is not compact but a surface as well as ones such that the metric completion is compact but not a surface.
For the first one consider the Euclidean plane~$\mathbb{R}^2$ and for the second one the baker's map surface studied by Chamanara in \cite{chamanara_04} or the icicled surface that is defined in \autoref{exa_icicled_surface}.

For a translation surface $(X,\mathcal{A})$ we call the points in $\overline{X} \setminus X$ \emph{singularities} of $(X,\mathcal{A})$.
In contrast to finite translation surfaces, very exotic behaviour of singularities can occur.
Consider for example the open disk $B(0,1) \subseteq \mathbb{R}^2$. It is a translation surface where the set of singularities is the unit circle, hence not discrete. This leads to the undesirable feature that for every point in the surface and every direction in $S^1$, the geodesic flow from that point in this direction is only defined for finite time.
This behaviour makes it quite hard to study dynamical properties, so we restrict ourselves to a certain class of translation surfaces in this article.

\begin{conv}[Translation surfaces have discrete singularities] \label{conv_discrete_singularities}
 For all translation surfaces in this article, the set of singularities is discrete in the metric completion, i.e.\ for each singularity $\sigma$, there exists an $\epsilon > 0$ such that $\sigma$ is the only singularity in $B(\sigma, \epsilon) \subseteq \overline{X}$. 
\end{conv}

\noindent
This convention makes it also possible to distinguish different types of singularities in the following way.

\begin{definition}[Cone angle, infinite angle, and wild singularities] \label{def_types_singularities}
 Let $(X,\mathcal{A})$ be a translation surface and $\sigma$ a singularity of $(X,\mathcal{A})$.
 \begin{enumerate}
  \item The singularity $\sigma$ is called \emph{cone angle singularity of multiplicity $k>0$} if there exist
  \begin{itemize}
   \item $\epsilon > 0$,
   \item an open neighbourhood $B$ of $\sigma$ in $\overline{X}$, and
   \item a $k$--cyclic translation covering from $B \setminus \{\sigma\}$ to the once-punctured disk ${B(0,\epsilon) \setminus \{0\} \subseteq \mathbb{R}^2}$.
  \end{itemize}
  
  If $k=1$ then the singularity $\sigma$ is also called \emph{removable singularity} or \emph{flat point}.

  \item The singularity $\sigma$ is called \emph{infinite angle singularity} or \emph{cone angle singularity of multiplicity $\infty$} if there exist
  \begin{itemize}
   \item $\epsilon > 0$,
   \item an open neighbourhood $B$ of $\sigma$ in $\overline{X}$, and
   \item an infinite cyclic translation covering from $B \setminus \{\sigma\}$ to the once-punctured disk $B(0,\epsilon) \setminus \{0\} \subseteq \mathbb{R}^2$.
  \end{itemize}

  \item The singularity $\sigma$ is called \emph{wild} if it is neither a cone angle nor an infinite angle singularity.
 \end{enumerate}
\end{definition}

\noindent
For cone angle and infinite angle singularities, all geometric information is encoded in the multiplicity.
To understand wild singularities in a similar way, we have to describe them in more detail. The first attempt at this was done by Bowman and Valdez in \cite{bowman_valdez_13}.
We recall it very briefly.

\begin{definition}[Space of linear approaches]
 Let $(X,\mathcal{A})$ be a translation surface, $x\in \overline{X}$, and $\epsilon > 0$. We define
 \begin{equation*}
  \mathcal{L}^\epsilon (x) \coloneqq \left\{ \gamma \colon (0,\epsilon) \to X : \gamma \text{ is a geodesic curve and } \lim\limits_{t \to 0} \gamma(t) = x \right\}
  .
 \end{equation*}
 
 If $x$ is a wild singularity then we can deduce from \autoref{def_types_singularities} that there exists no $\epsilon > 0$ such that all geodesic curves starting in $x$ can be extended to have at least length~$\epsilon$.
 Therefore, we consider equivalence classes instead of curves:
 $\gamma_1 \in \mathcal{L}^\epsilon (x)$ and $\gamma_2 \in \mathcal{L}^{\epsilon'} (x)$ are called \emph{equivalent} if $\gamma_1 (t) = \gamma_2 (t)$ for every $t\in (0,\min\{\epsilon, \epsilon' \})$.

 The space
 \begin{equation*}
  \mathcal{L} (x) \coloneqq \bigsqcup\limits_{\epsilon > 0} \mathcal{L}^\epsilon (x) / \sim
 \end{equation*}
 is called the \emph{space of linear approaches of $x$} and the equivalence class $[\gamma]$ of $\gamma\in \mathcal{L}^{\epsilon}(x)$ is called a \emph{linear approach} to the point $x$.
\end{definition}

\noindent
Every $\mathcal{L}^{\epsilon}(x)$ can be embedded in $\mathcal{L}(x)$ and also in $\mathcal{L}^{\epsilon'}(x)$ for every $\epsilon'>0$ with $\epsilon' \leq \epsilon$.
The family of all spaces $\mathcal{L}^{\epsilon}(x)$ together with these embeddings is a direct system and $\mathcal{L}(x)$ is the colimit of this direct system.

We now describe the topology with which we endow the space of linear approaches.

\begin{definition}[Topology on $\mathcal{L}(x)$]
 Let $(X, \mathcal{A})$ be a translation surface and $x \in \overline{X}$.
 For $\epsilon > 0$ we can define the \emph{uniform metric} $d_\epsilon$ on $\mathcal{L}^\epsilon(x)$ using the translation metric $d_X$ on $X$:
 \begin{equation*}
  d_\epsilon (\gamma_1, \gamma_2) = \sup\limits_{0 < t < \epsilon} d_X \left( \gamma_1 (t), \gamma_2(t) \right) .
 \end{equation*}
 The uniform metric defines a topology on $\mathcal{L}^\epsilon (x)$.
 As $\mathcal{L}(x)$ is the colimit of the direct system of the spaces $\mathcal{L}^\epsilon(x)$, we can define the \emph{final topology} on $\mathcal{L}(X)$, i.e.\ the finest topology such that all embeddings $\mathcal{L}^\epsilon(X) \hookrightarrow \mathcal{L}(X)$ are continuous.
\end{definition}

\noindent
The last concept we recall from \cite{bowman_valdez_13} is is that of a rotational component as a class of linear approaches with a one-dimensional translation structure on it.
The idea is to mimic the notion of ``going around a singularity'' for cone angle and infinite angle singularities.

\begin{figure}
\begin{center}
\begin{tikzpicture}
 \draw (0,0) node[left,align=center] {base\\point} -- (-25:2) node[below] {$\epsilon \cdot e^{a\cdot \mathrm{i}}$} arc (-25:75:2) node[above] {$\epsilon \cdot e^{b\cdot \mathrm{i}}$} -- (0,0);
 \foreach \x in {-15,-5,...,65}
  \draw[very thin] (0,0) -- (\x:2); 

 \begin{scope}[xshift=5.1cm]
  \draw (0,0) -- (-55:2) arc (-55:195:2) -- (0,0);
  \foreach \x in {-45,-35,...,185}
   \draw[very thin] (0,0) -- (\x:2); 
 \end{scope}

 \begin{scope}[xshift=10.2cm]
  \draw (-35:2) arc (-35:325:2);
  \draw (325:2) arc (-35:51:2) -- (0,0);
  \draw[densely dotted] (0,0) -- (-35:2);
  \foreach \x in {-25,-14.859,...,46}
   \draw[densely dotted, very thin] (0,0) -- (\x:2);
  \foreach \x in {56.141,66.282,...,410}
   \draw[very thin] (0,0) -- (\x:2); 
 \end{scope}
\end{tikzpicture}
\end{center}
\label{fig_examples_angular_sectors}
\caption{Examples for $i_\epsilon (U)$ with $I=[a,b]$: in the first two examples we have $b-a < 2 \pi$ whereas in the last example we have $b-a> 2\pi$.}
\end{figure}
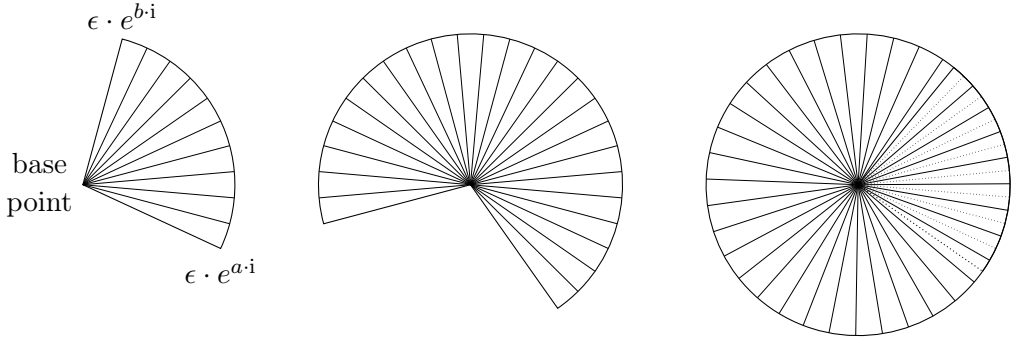

For $\epsilon > 0$ and a \emph{generalized interval} $I$, i.e.\ a nonempty connected subset of $\mathbb{R}$, we consider the infinite strip $\left\{ z \in \mathbb{C} : \Re(z) < \log \epsilon, \Im(z) \in I \right\} \subseteq \mathbb{C}$. 
Via the injective map 
\begin{equation*}
 f \colon \mathbb{C} \to \left( \mathbb{C}\setminus \{0\} \right) \times \mathbb{R}, \ z \mapsto \left( e^z, \Im(z) \right)
 ,
\end{equation*}
we can spiral the strip around the puncture at $0$. We endow the image $U$ of the strip under $f$ with the pullback of the Euclidean metric on $\mathbb{C}$ via the projection to the first component.
Some examples of the projection of~$U$ are sketched in \autoref{fig_examples_angular_sectors}.

Such an image $U$ together with an embedding in $X$ defines an angular sector, as is made precise in the following definition.

\begin{definition}[Angular sector]
 Let $(X, \mathcal{A})$ be a translation surface. An \emph{angular sector} is a triple $(I, \epsilon, i_\epsilon)$ of a generalized interval $I$, an $\epsilon > 0$, and an isometric embedding $i_\epsilon$ of 
 \begin{equation*}
  U \coloneqq f \left( \left\{ z \in \mathbb{C} : \Re(z) < \log \epsilon, \Im(z) \in I \right\} \right)
 \end{equation*}
 into $X$.
\end{definition}

\noindent
For an angular sector $(I,\epsilon, i_\epsilon)$ and $y\in I$, the limit
  $\lim_{x \to -\infty} (i_\epsilon \circ f)(x+\mathrm{i}y)$
is a point in~$\overline{X}$ and independent of $y$.
This point is called the \emph{basepoint} of the angular sector $(I,\epsilon, i_\epsilon)$.

If $(I,\epsilon,i_\epsilon)$ is an angular sector with basepoint $x \in \overline{X}$
then we can define a map $f_{(I,\epsilon,i_\epsilon)} \colon I\rightarrow \mathcal{L}(x)$ as follows.
For every $y \in I$, the image of the map 
\begin{equation*}
 i_\epsilon \circ f \colon \left\{ z \in \mathbb{C} : \Re(z) < \log \epsilon, \Im(z)=y \right\} \to X
\end{equation*}
is a geodesic segment of length $\epsilon$ and hence induces an element in $\mathcal{L}^\epsilon(x)$. Let $f_{(I,\epsilon,i_\epsilon)} (y)$ be the corresponding linear approach in $\mathcal{L}(x)$.

In the next definition, we use angular sectors to define an equivalence relation on the space $\mathcal{L}(x)$ of linear approaches of a point~$x$. In an informal way, we can describe it as two linear approaches being ``contained'' in the image of the same angular sector.

\begin{definition}[Rotational component] \label{def_rotational_component}
 Let $(X, \mathcal{A})$ be a translation surface, $x\in \overline{X}$, and $[\gamma_1], [\gamma_2] \in \mathcal{L}(x)$.
 The linear approaches $[\gamma_1]$ and $[\gamma_2]$ are called \emph{$R$--equi\-va\-lent} if there exists an angular sector $(I, \epsilon, i_\epsilon)$ with basepoint $x$ and $y_1, y_2 \in I$ such that $f_{(I,\epsilon,i_\epsilon)} (y_1)=[\gamma_1]$ and $f_{(I,\epsilon,i_\epsilon)} (y_2) = [\gamma_2]$.

 The $R$--equivalence class $\overline{[\gamma]}$ of $[\gamma] \in \mathcal{L}(x)$ is called a \emph{rotational component} of~$x$.
\end{definition}

\noindent
We say that a linear approach $[\gamma]$ is \emph{contained} in a rotational component~$c$ if $\overline{[\gamma]} = c$ holds.
The set of all linear approaches assigned to a given angular sector $(I, \epsilon, i_\epsilon)$ can be described as
\begin{equation*}
 V(I, \epsilon, i_\epsilon) \coloneqq \left\{ f_{(I,\epsilon,i_\epsilon)} (y) : y\in I \right\} .
\end{equation*}
Then $V(I, \epsilon, i_\epsilon)$ is a subset of the space $\mathcal{L}(x)$ of linear approaches. All linear approaches in $V(I, \epsilon, i_\epsilon)$ are contained in the same rotational component~$c$, where $c \coloneqq \overline{f_{(I,\epsilon,i_\epsilon)} (y)}$ for an arbitrary $y\in I$.
In other words, the map 
\begin{equation*}
 \varphi_{(I, \epsilon, i_\epsilon)} \colon V(I, \epsilon, i_\epsilon) \to \mathbb{R}, \ f_{(I,\epsilon,i_\epsilon)} (y) \mapsto y
\end{equation*}
is inverse to $f_{(I,\epsilon,i_\epsilon)}$ and hence $V(I,\epsilon, i_\epsilon)$ and $I$ are in one-to-one correspondence.

We want to define a topology on a rotational component $\overline{[\gamma]}$ of $x$ by considering all angular sectors $(I, \epsilon, i_\epsilon)$ with basepoint $x$ such that $I$ is an open interval and $f_{(I,\epsilon,i_\epsilon)} (y)$ is contained in $\overline{[\gamma]}$ for every $y \in I$.
The union of the images of all such $f_{(I, \epsilon, i_\epsilon)}$ covers the rotational component except for possibly two linear approaches $[\gamma_{left}]$ and $[\gamma_{right}]$.
In this case, half-closed intervals $I = (a,b]$ with $f_{(I, \epsilon, i_\epsilon)} (b) \in \{ [\gamma_{left}], [\gamma_{right}] \}$ have to be allowed as open sets to obtain a cover of $\overline{[\gamma]}$.

We choose the collection of the corresponding sets $V(I,\epsilon, i_\epsilon)$ as a basis of the topology. In particular, any such $\varphi_{(I, \epsilon, i_\epsilon)}$ for an open interval $I$ is a homeomorphism.

\begin{definition}[Translation structure on rotational components] \label{def_translation_structure_rotational_component}
 Let $(X, \mathcal{A})$ be a translation surface, $x \in \overline{X}$, and $\overline{[\gamma]}$ a rotational component of $x$ that contains more than one linear approach. Then
 \begin{align*}
  \{ (V(I, \epsilon, i_\epsilon), \varphi_{(I, \epsilon, i_\epsilon)}):\
  & (I,\epsilon, i_\epsilon) \text{ an angular sector with basepoint } x, \\
  &  I \text{ open},\ f_{(I, \epsilon, i_\epsilon)} (I) \subseteq \overline{[\gamma]} \} \\
  \cup\ \{ (V(I, \epsilon, i_\epsilon), \varphi_{(I, \epsilon, i_\epsilon)}):\
  & (I,\epsilon, i_\epsilon) \text{ an angular sector with basepoint } x, \\
  &  I = (a,b], \ f_{(I, \epsilon, i_\epsilon)} (I) \subseteq \overline{[\gamma]}, \\
  &  f_{(I, \epsilon, i_\epsilon)}  (b) \in \{ [\gamma_{left}], [\gamma_{right}] \} \}
 \end{align*}
 forms an atlas of $\overline{[\gamma]}$.

 Therefore, $\overline{[\gamma]}$ is a one-dimensional manifold, possibly with boundary. As the transition functions are actually translations in $\mathbb{R}$, $\overline{[\gamma]}$ even carries a one-dimensional translation structure.
\end{definition}

\noindent
Using the one-dimensional translation structure on $\overline{[\gamma]}$, we can pull back the Euclidean metric from $\mathbb{R}$ to $\overline{[\gamma]}$ and for this we can measure the \emph{length of a rotational component}.

\section{Criterion for infinite genus} \label{sec_criterion_infinite_genus}

The purpose of this article is to show that a translation surface with certain conditions has infinite genus. We start by giving a feasible criterion for the infinity of genus in this section.

The genus of a surface is often defined by Betti numbers or by the idea of handles of a surface. Instead of those approaches, we will employ a definition in this section that uses the number of nonseparating curves.
For a connected surface $X$, a simple closed curve $\gamma \colon [0,l] \to X$ is called \emph{nonseparating} if $X$ with the image of $\gamma$ removed is connected. In the same way, we say that a set of simple closed curves $\{\gamma_1, \ldots, \gamma_n\}$ is \emph{nonseparating} if~$X$ with the images of $\gamma_1,\ldots,\gamma_n$ removed is connected.

Similarly, we say that a closed curve $\gamma$ in $\overline{X}$ is \emph{nonseparating} if $\overline{X} \setminus \im(\gamma)$ is connected and a set of closed curves $\{\gamma_1, \ldots, \gamma_n\}$ in $\overline{X}$ is \emph{nonseparating} if $\overline{X} \setminus \im(\gamma_1) \cup \ldots \cup \im(\gamma_m)$ is connected.
Note that the notion of a nonseparating set of curves is in both cases stronger than being a set of nonseparating curves.

\begin{definition}[Genus]\label{def_genus}
 Let $X$ be a connected, orientable surface. The \emph{genus} of $X$ is $g \in \mathbb{N}$ if the following equivalent conditions are true.
 \begin{enumerate}
  \item The maximum cardinality of a nonseparating set of disjoint curves in $X$ is $g$.
  \item The maximum cardinality of a nonseparating set of curves in $X$ is~$2g$.
 \end{enumerate}
\end{definition}

\noindent
Within the meaning of the previous definition, a connected, orientable surface $X$ is said to have \emph{infinite genus} if for every $n \in \mathbb{N} = \{0,1, \ldots\}$ there exists a nonseparating set of curves in $X$ with cardinality $n$. This is equivalent to $X$ containing subsurfaces of arbitrarily large genus.

In \autoref{prop_criterion_saddle_connections_infinite_genus}, we specialize this definition to a criterion for the case of translation surfaces using \emph{saddle connections}, i.e.\ geodesic curves $\gamma \colon [0,l] \to \overline{X}$ such that $\gamma((0,l))$ is contained in $X$ and $\gamma(0), \gamma(l) \in \overline{X} \setminus X$.
The criterion involves so-called \emph{left-to-right curves} that connect one side of a given curve to the other side. This notion is made precise in \autoref{def_lefttoright_curves_curves} and \autoref{def_lefttoright_curves_saddle_connections}.

\begin{prop}[Saddle connections and infinite genus] \label{prop_criterion_saddle_connections_infinite_genus}
 Let $(X,\mathcal{A})$ be a translation surface and $\sigma$ a singularity. Suppose that for every $n \geq 1$ there exists a set of $n$ saddle connections from $\sigma$ to itself such that the saddle connections intersect exactly in $\sigma$ and the set has left-to-right curves. Then $X$ has infinite genus.
\end{prop}

\noindent
We will prove the proposition by several lemmas and start with a precise definition of left-to-right curves.

\begin{definition}[Left-to-right curves of curves in $\boldsymbol{X}$] \label{def_lefttoright_curves_curves}
 Let $X$ be a connected, orientable surface, $n\geq 1$, and $\gamma_1,\ldots, \gamma_n$ simple closed curves in~$X$ that intersect pairwise in exactly one point $x\in X$.
 \begin{enumerate}
  \item Let $\epsilon > 0$ be small enough such that the $\epsilon$--neighbourhood $N$ of $\gamma_1$ is a tubular neighbourhood.
  Then $N$ is topologically an annulus. So $N \setminus \im(\gamma_1)$ consists of two connected components $N_l$ and $N_r$. Considering the underlying orientation of $\im(\gamma_1)$ we call points in $N_l$ and $N_r$ \emph{points on the left of $\gamma_1$} and \emph{points on the right of $\gamma_1$}, respectively.

  \item A curve in $X \setminus \im(\gamma_1)$ from a point on the left of $\gamma_1$ to a point on the right of $\gamma_1$ is called a \emph{left-to-right curve} of $\gamma_1$.

  \item Choose a tubular neighbourhood $N$ of $\gamma_n$ and let $N_l$ and $N_r$ be as before.
  We have that $N_l \setminus \left(\im(\gamma_1)\cup \ldots \cup \im(\gamma_{n-1})\right)$ consists of one or more connected components (see \autoref{fig_crossing_nonseparating_curves}).
  The boundary of such a connected component consists of a subset of the boundary of~$N$ and of subsets of the images of some $\gamma_i$.
  As the curves $\gamma_1,\ldots, \gamma_n$ intersect in exactly one point, there is only one connected component $N_l^\ast$ whose boundary contains $\im(\gamma_n)$.

  We call a point in this connected component $N_l^\ast$ a \emph{point on the left of $\gamma_n$ with respect to $\gamma_1,\ldots, \gamma_{n-1}$}. Similarly, we define $N_r^\ast$ and \emph{points on the right of $\gamma_n$ with respect to $\gamma_1,\ldots, \gamma_{n-1}$}.
  Then a curve in $X \setminus ( \im(\gamma_1) \cup \ldots \cup \im(\gamma_n))$ is called a \emph{left-to-right curve of $\gamma_n$ with respect to $\gamma_1,\ldots, \gamma_{n-1}$} if it connects a point on the left of $\gamma_n$ with respect to $\gamma_1,\ldots, \gamma_{n-1}$ to a point on the right of $\gamma_n$ with respect to $\gamma_1,\ldots, \gamma_{n-1}$.

  \begin{figure}
  \begin{center}
  \begin{tikzpicture}
   \draw (0,0) to[bend left=20] (4,6) node[above] {$\gamma_n$};
   \draw[color=gray] (-0.5,0.2) to[bend left=20] (3.6,6.2);
   \draw[color=gray] (0.5,-0.2) to[bend left=20] node[right,pos=0.9] {$N$} (4.4,5.8);

   \draw (-1,1.73) node[left] {$\gamma_{n-1}$} to[bend left=10] (5,5);
   \draw (-1,2.9) node[left] {$\vdots$};
   \draw (-1,4) node[left] {$\gamma_2$} to[bend right=5] (5,3);
   \draw (-0.6,5.8) node[left] {$\gamma_1$} to[bend left=5] (3.13,1);
  \end{tikzpicture}
  
  \label{fig_crossing_nonseparating_curves}
  \caption{The connected components $N_l^\ast$ and $N_r^\ast$ are the ones that are bounded by $\gamma_n$ and segments of the gray curve and of some $\gamma_1,\ldots,\gamma_{n-1}$}
  \end{center}
\end{figure}

  \item We say that the set of curves $\{\gamma_1,\ldots, \gamma_n\}$ \emph{has left-to-right curves} if every curve has a left-to-right curve with respect to the other ones.
 \end{enumerate}
\end{definition}

\noindent
For a translation surface $(X,\mathcal{A})$ and a curve in $\overline{X}$, for instance a saddle connection, we cannot find such a tubular neighbourhood $N$ as in \autoref{def_lefttoright_curves_curves} (i) but use a slightly different neighbourhood instead. So we can define left-to-right curves for a special type of curves in $\overline{X}$ in a similar way while avoiding the singularity.

\begin{definition}[Left-to-right curves of curves in $\boldsymbol{\overline{X}}$] \label{def_lefttoright_curves_saddle_connections}
 Let $(X, \mathcal{A})$ be a translation surface, $\sigma$ a singularity, $n\geq 1$, and $\gamma_1,\ldots, \gamma_n$ simple closed curves in $X \cup \{\sigma\}$ whose images contain $\sigma$ exactly as startpoint and endpoint and that are disjoint in their interiors.

 \begin{enumerate}
  \item Let $l$ be the length of $\gamma_1$ and let $\epsilon > \epsilon'>0$ be sufficiently small. Consider the set $N \subseteq X \cup \{\sigma\}$ which is the union of $B(\sigma, \epsilon)$ and the open $\epsilon'$--neighbourhood $\widetilde{N}$ of the segment $\gamma_1([\epsilon, l-\epsilon])$. Again, $\widetilde{N} \setminus \im(\gamma_1)$ consists of two connected components $\widetilde{N}_l$ and $\widetilde{N}_r$. In this situation, we call points in $\widetilde{N}_l$ \emph{points on the left of $\gamma_1$} and points in~$\widetilde{N}_r$ \emph{points on the right of $\gamma_1$}, with respect to the orientation of~$\gamma_1$.

  \item A curve in $X\setminus \im(\gamma_1)$ from a point on the left of $\gamma_1$ to a point on the right of $\gamma_1$ is called a \emph{left-to-right curve} of $\gamma_1$.

  \item Choose a neighbourhood $N$ of $\gamma_n$ and let $\widetilde{N}_l$ and $\widetilde{N}_r$ be as before.
  Again, it is possible that $\widetilde{N}_l \setminus \left(\im(\gamma_1)\cup \ldots \cup \im(\gamma_{n-1})\right)$ consists of more than one connected component. In this case, we can avoid the indicated behaviour by choosing $\epsilon'>0$ small enough such that none of the curves $\gamma_1, \ldots, \gamma_{n-1}$ intersects $\widetilde{N}_l$ and $\widetilde{N}_r$.
  We call a point in these newly chosen $\widetilde{N}_l$ and $\widetilde{N}_r$ a \emph{point on the left of $\gamma_n$ with respect to $\gamma_1,\ldots, \gamma_{n-1}$} and a \emph{point on the right of $\gamma_n$ with respect to $\gamma_1,\ldots, \gamma_{n-1}$}, respectively.

  Then a curve in $X \setminus ( \im(\gamma_1) \cup \ldots \cup \im(\gamma_n))$ is called a \emph{left-to-right curve of $\gamma_n$ with respect to $\gamma_1,\ldots, \gamma_{n-1}$} if it connects a point on the left of $\gamma_n$ with respect to $\gamma_1,\ldots, \gamma_{n-1}$ to a point on the right of $\gamma_n$ with respect to $\gamma_1,\ldots, \gamma_{n-1}$.

  \item We say that the set of curves $\{\gamma_1,\ldots, \gamma_n\}$ \emph{has left-to-right curves} if every curve has a left-to-right curve with respect to the other ones.
 \end{enumerate}
\end{definition}

\noindent
Note that the existence of left-to-right curves as in \autoref{def_lefttoright_curves_curves} and \autoref{def_lefttoright_curves_saddle_connections} does not depend on the choice of $\epsilon$ or $\epsilon'$ as long as these values are small.

By means of left-to-right curves we can now formulate a criterion for a set of curves to be nonseparating.

\begin{lem}[Criterion for simple closed curves to be nonseparating] \label{lem_criterion_simple_closed_curves_nonseparating}
 Let~$X$ be a connected, orientable surface, $n\geq 1$, and $\gamma_1,\ldots, \gamma_n$ simple closed curves in $X$ that intersect pairwise in exactly one point $x\in X$.
 Then the set $\{\gamma_1,\ldots,\gamma_n\}$ is nonseparating if and only if the set has left-to-right curves.

\begin{proof}
 We prove this statement by induction on the number $n$ of curves. 
 For the base case we show that the curve $\gamma_1$ is nonseparating if and only if it has a left-to-right curve.
 
 If $\gamma_1$ is nonseparating then we can take any point on the left and any point on the right of $\gamma_1$ and, as $X \setminus \im(\gamma_1)$ is connected, there exists a curve connecting these two points without intersecting $\gamma_1$.

 Now assume we have such a left-to-right curve $\delta_1$ connecting a point $x_l$ on the left of~$\gamma_1$ and a point $x_r$ on the right of $\gamma_1$. Choose two points $x_1, x_2 \in X \setminus \im(\gamma_1)$. We have to show that there exists a curve $\beta$ in $X\setminus \im(\gamma_1)$ that connects $x_1$ and $x_2$.
  
 Let $N$, $N_l$ and $N_r$ be as in \autoref{def_lefttoright_curves_curves}. As $X$ is connected, there exists a curve $\beta'$ in~$X$ that connects $x_1$ to $x_2$. If $\beta'$ does not intersect $\gamma_1$ then we can choose $\beta \coloneqq \beta'$.
 If it does then let $\beta_+'$ be the subcurve of $\beta'$ from $x_1$ to the first intersection of $\beta'$ and $\overline{N}$ and let $\beta_-'$ be the subcurve of $\beta'$ from the last intersection of $\beta'$ and $\overline{N}$ to $x_2$ (see \autoref{fig_left_and_right_of_curve} for a sketch).

 \begin{figure}
  \begin{center}
  \begin{tikzpicture}
   \draw (0,0) to[bend right=20] (4,6) node[above] {$\gamma_1$};
   \draw[color=gray, name path=leftboundary] (-0.4,0.4) to[bend right=20] (3.6,6.2);
   \draw[color=gray, name path=rightboundary] (0.4,-0.4) to[bend right=20] node[right,pos=0.8] {$N$} (4.4,5.8);

   \fill (-0.5,3) node[left] {$x_1$} circle (1pt);
   \fill (3.5,1) node[right] {$x_2$} circle (1pt);
   \draw[name path=betaprime] (-0.5,3) .. controls +(60:0.8cm) and +(-180:0.5cm) ..
   (1.5,4.5) node[above] {$\beta'$} .. controls +(0:0.8cm) and +(80:1.5cm) .. node[left] {$\beta'_+$}
   (1.3,2) .. controls +(-100:3.2cm) and +(95:1.7cm) ..
   (2.4,1.3) .. controls +(-85:1.5cm) and +(-130:1.5cm) .. node[below] {$\beta'_-$}
   (3.5,1);
   \fill[name intersections={of=leftboundary and betaprime, by={endpointleft}}] (endpointleft) circle (1pt);
   \fill[name intersections={of=rightboundary and betaprime, by={startpointright}}] (startpointright) circle (1pt);
  \end{tikzpicture}
  
  \label{fig_left_and_right_of_curve}
  \caption{This configuration of $\gamma_1$ and $\beta'$ is treated in Case 2.}
  \end{center}
 \end{figure}
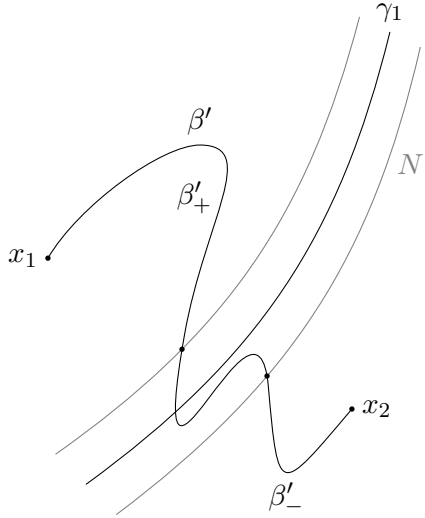
  
 Case 1: The endpoint of $\beta_+'$ and the startpoint of $\beta_-'$ both belong to $\overline{ N_r}$ or both to~$\overline{ N_l}$. Then we can choose a curve between these two points in the connected set $N_r$ or $N_l$ and the concatenation of this curve with $\beta_+'$ and $\beta_-'$ gives us a curve $\beta$ as desired.
  
 Case 2: The endpoint of $\beta_+'$ belongs to $\overline{ N_l}$ or to $\overline{ N_r}$ and the startpoint of $\beta_-'$ belongs to the other.
 For simplicity of notation, let $\beta_+'$ end in $\overline{ N_l}$ and $\beta_-'$ start in $\overline{ N_r}$. Again, as $N_l$ and $N_r$ are connected we find curves connecting the endpoint of $\beta_+'$ to~$x_l$ in $N_l$ and connecting $x_r$ to the startpoint of $\beta_-'$ in $N_r$. Concatenating all these curves with the left-to-right curve $\delta_1$ in the correct order yields a curve $\beta$ as desired.
 
 This concludes the proof of the base case. For the inductive step, let $n\geq 2$ and $\{\gamma_1,\ldots, \gamma_n\}$ be a set of simple closed curves that intersect pairwise exactly in $x$ and are such that the set $\{\gamma_1,\ldots, \gamma_{n-1}\}$ is nonseparating if and only if the set has left-to-right curves.

 Let $\{\gamma_1,\ldots, \gamma_n\}$ be a nonseparating set of curves.
 We show for every $i\in \{1,\ldots,n\}$ that $\gamma_i$ has a left-to-right curve with respect to $\gamma_1,\ldots,\gamma_{i-1},\gamma_{i+1}, \ldots, \gamma_n$ by choosing a point on the left and a point on the right of~$\gamma_i$ with respect to $\gamma_1,\ldots,\gamma_{i-1}, \gamma_{i+1}, \ldots, \gamma_n$. As $X \setminus \{\im(\gamma_1)\cup \ldots \cup \im(\gamma_n)\}$ is connected, there exists a curve connecting the two chosen points without intersecting one of the curves $\gamma_1,\ldots, \gamma_n$.

 Now assume that the set $\{\gamma_1,\ldots,\gamma_n\}$ has left-to-right curves.
 Then $\{\gamma_1,\ldots, \gamma_{n-1}\}$ has left-to-right curves and is hence nonseparating. Let $\delta_n$ be the left-to-right curve of $\gamma_n$ with respect to $\gamma_1,\ldots, \gamma_{n-1}$ connecting a point $x_l$ on the left of $\gamma_n$ to a point $x_r$ on the right of $\gamma_n$.

 We have to show that for two chosen points $x_1, x_2 \in X \setminus \left(\im(\gamma_1)\cup \ldots \cup \im(\gamma_n)\right)$ there exists a curve $\beta_n$ in $X \setminus \left(\im(\gamma_1)\cup \ldots \cup \im(\gamma_n)\right)$ that connects $x_1$ and $x_2$.
 As the set $\{\gamma_1, \ldots, \gamma_{n-1}\}$ is nonseparating, there exists a curve $\beta_n'$ in $X \setminus \left(\im(\gamma_1) \cup \ldots \cup \im(\gamma_{n-1})\right)$ that connects $x_1$ and $x_2$.

 If $\beta_n'$ does not intersect $\gamma_n$ then we can choose $\beta_n \coloneqq \beta_n'$. If it does then it also intersects one of the connected components $N_l^\ast$ or $N_r^\ast$ of $N \setminus \left(\im(\gamma_1)\cup \ldots \cup \im(\gamma_n)\right)$. This is because the boundaries of all other connected components contain only one point of $\im(\gamma_n)$ which is the point $x= \im(\gamma_1) \cap \ldots \cap \im(\gamma_n)$.
 Now let $\beta_+'$ be the subcurve of $\beta_n'$ from $x_1$ to the first intersection of $\beta_n'$ and $\overline{ N_l^\ast \cup N_r^\ast }$ and let $\beta_-'$ be the subcurve of $\beta_n'$ from the last intersection of $\beta_n'$ and $\overline{ N_l^\ast \cup N_r^\ast }$ to $x_2$.
 Then we proceed as in the base case and construct a curve $\beta_n$ in $X \setminus \left(\im(\gamma_1)\cup \ldots \cup \im(\gamma_n)\right)$ that connects $x_1$ and $x_2$.

 This concludes the proof of the inductive step and hence the proof of the lemma.
\end{proof}
\end{lem}

\noindent
It follows by the same arguments that the criterion in \autoref{lem_criterion_simple_closed_curves_nonseparating} is also true for saddle connections or, more generally, simple closed curves in $\overline{X}$ whose image contains singularities exactly as startpoint and endpoint. While the beginning of the proof of the base case is literally the same, we have to choose an $\epsilon > 0$ small enough such that all intersection points of $\beta'$ and~$\gamma_1$ are contained in the $\epsilon'$--neighbourhood $\widetilde{N}$ of the segment $\gamma_1([\epsilon, l-\epsilon])$. Then we can use $\widetilde{N}$ instead of $N$ and finish the proof of the base case in the same way as before. The replacement of $N$ by $\widetilde{N}$ also makes the proof of the inductive step work for this type of curves as the connected components of $\widetilde{N} \setminus \im(\gamma_1)$ can also play the role of the connected components $N_l^\ast$ and~$N_r^\ast$.

Altogether, we have the following version of \autoref{lem_criterion_simple_closed_curves_nonseparating}.

\begin{lem}[Criterion for curves in $\overline{X}$ to be nonseparating] \label{lem_criterion_saddle_connections_nonseparating}
 Let $(X, \mathcal{A})$~be a translation surface, $\sigma$ a singularity, $n\geq 1$, and $\gamma_1,\ldots, \gamma_n$ simple closed curves in $X \cup \{\sigma\}$ whose images contain $\sigma$ exactly as startpoint and endpoint and that are disjoint in their interiors.
 Then the set $\{\gamma_1,\ldots,\gamma_n\}$ is nonseparating if and only if the set has left-to-right curves.
\end{lem}

\noindent
The statement of \autoref{lem_criterion_saddle_connections_nonseparating} is already close to our goal.
However, as ``genus'' is a concept for surfaces we have to consider curves in $X$ instead of curves in $\overline{X}$ to determine infinite genus. Because of this, we show in \autoref{lem_nonseparating_curves_in_metric_completion} how to replace curves in $\overline{X}$ by curves in $X$ without disturbing their left-to-right curves.

In the proof of the lemma, we use some curves whose existence is ensured by the following proposition.

\begin{prop}[Balls around singularities are path-connected]\label{prop_balls_around_singularities_pathconnected}
 Let $(X,\mathcal{A})$ be a translation surface, $\sigma$ a singularity, and $\epsilon > 0$. Then $B(\sigma,\epsilon)\cap X$ is path-connected.

\begin{proof}
 For a fixed $\epsilon > 0$, let $x,y\in B(\sigma, \epsilon) \cap X$ and $\delta>0$ small enough such that $d(x,\sigma) < \epsilon - 3\delta$ and $d(y,\sigma) < \epsilon - 3\delta$.
 As $\sigma$ is contained in the metric completion of $X$, there exists a point $z\in X$ with $d(\sigma, z) < \delta$.
 Then we have
 \begin{equation*}
  d(x,z) \leq d(x,\sigma) + d(\sigma, z) < \epsilon - 2\delta
 \end{equation*}
 and $d(y,z) < \epsilon - 2\delta$.

 As the metric of $X$ is the path-length metric, there exists a curve~$\gamma_x$ in $X$ which connects $x$ to $z$ and has at most length $\epsilon - \delta$. This means that, for every point in the image of $\gamma_x$, the distance to $\sigma$ is at most $d(\sigma, z) + (\epsilon - \delta) < \epsilon$. Hence, $\gamma_x$ is a curve in $B(\sigma, \epsilon)\cap X$.

 In the same way we can define a curve $\gamma_y$ in $B(\sigma, \epsilon) \cap X$ from $y$ to~$z$ and, by concatenation of $\gamma_x$ with the reversed curve of $\gamma_y$, we obtain a curve from $x$ to $y$ in $B(\sigma, \epsilon) \cap X$. This means that $B(\sigma, \epsilon) \cap X$ is path-connected.
\end{proof}
\end{prop}

\noindent
The path-connectedness is also fulfilled for $B(\sigma, \epsilon) \subseteq \overline{X}$. The proof is literally the same as the metric of $\overline{X}$ is also the path-length metric.

\begin{lem}[Nonseparating curves in $\overline{X}$ give rise to nonseparating curves in $X$] \label{lem_nonseparating_curves_in_metric_completion}
 Let $(X,\mathcal{A})$ be a translation surface, $\sigma$ a singularity, $n\geq 1$, and $\{ \gamma_1,\ldots,\gamma_n \}$ a nonseparating set of simple closed curves in $X \cup \{\sigma\}$ whose images contain $\sigma$ exactly as startpoint and endpoint and that are disjoint in their interiors.
 Then there also exists a set of simple closed curves $\{ \gamma_1',\ldots,\gamma_n' \}$ in $X$ that is nonseparating.

\begin{proof}
 By \autoref{lem_criterion_saddle_connections_nonseparating}, there exist left-to-right curves~$\delta_i$ of $\gamma_i$ with respect to $\gamma_1,\ldots, \gamma_{i-1}, \allowbreak \gamma_{i+1}, \ldots, \gamma_n$ for every $i\in \{1,\ldots,n\}$. 
 Choose $\epsilon > 0$ small enough such that the $\epsilon$--neigh\-bor\-hood $B(\sigma,\epsilon)$ of $\sigma$ avoids all $\delta_i$ and such that $\partial B(\sigma, \epsilon)$ intersects $\im(\gamma_i)$ at least two times for every $i\in \{1,\ldots, n\}$. For every $i \in \{1,\ldots,n\}$, the first intersection point of $\partial B(\sigma, \epsilon)$ and $\im(\gamma_i)$ (with respect to the orientation of $\im(\gamma_i)$) is called~$x_i^+$ and the last intersection point is called $x_i^-$.
 
 We now replace the curves $\gamma_1,\ldots,\gamma_n$ in $X \cup \{\sigma\}$ by curves in $X$ that have similar properties.
 For this choose a point $z\in B(\sigma, \epsilon) \setminus \{\sigma\}$ that will play the role of the current intersection point~$\sigma$.
 Because of the path-connectedness of $B(\sigma, \epsilon) \cap X$ (see \autoref{prop_balls_around_singularities_pathconnected}) we have a curve in $B(\sigma,\epsilon)\setminus \{\sigma\}$ from~$z$ to $x_1^+$ and a curve from $x_1^-$ to $z$ (see \autoref{fig_pushing_the_curves_out_of_the_singularity}).
 Now let~$\gamma_1'$ be the closed curve that is the concatenation of the curve from~$z$ to $x_1^+$, the subcurve of~$\gamma_1$ from $x_1^+$ to $x_1^-$, and the curve from~$x_1^-$ to~$z$.
 If~$\gamma_1'$ intersects itself then we smooth the crossing by joining other pairs of subcurves at the crossing. Hence, we can assume that~$\gamma_1'$ is simple.
 Also, as $\epsilon$ is chosen small enough, the curve $\delta_1$ is still a left-to-right curve of $\gamma_1'$ with respect to $\gamma_2, \ldots, \gamma_n$.

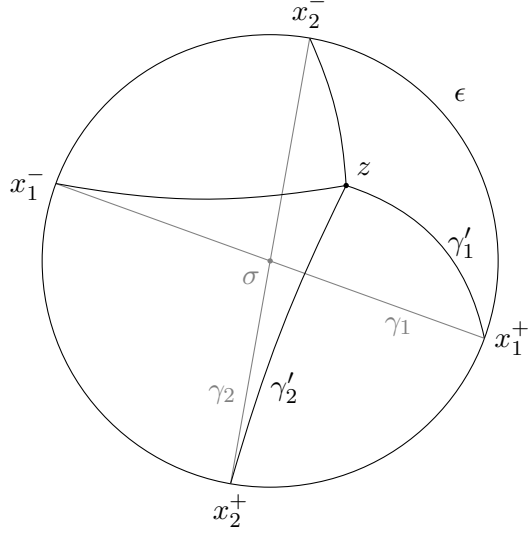
\begin{figure}
\begin{center}
 \begin{tikzpicture}
  \fill[color=gray] (0,0) node[below left] {$\sigma$} circle (1pt);
  \draw (0,0) circle (3cm);
  \draw[color=gray] (160:3) -- node[below, pos=0.8] {$\gamma_1$} (-20:3);
  \draw[color=gray] (80:3) -- node[left, pos=0.8] {$\gamma_2$} (-100:3);
  \fill (1,1) node[above right] {$z$} circle (1pt);
  \draw (1,1) to[bend left=10] (160:3) node[left] {$x_1^-$};
  \draw (1,1) to[bend left] node[above=0.1cm, pos=0.7] {$\gamma_1'$} (-20:3) node[right] {$x_1^+$};
  \draw (1,1) to[bend right=10] (80:3) node[above] {$x_2^-$};
  \draw (1,1) to[bend right=5] node[right, pos=0.7] {$\gamma_2'$} (-100:3) node[below] {$x_2^+$};
  \draw (2.5,2.2) node {$\epsilon$};
  \end{tikzpicture}
 \label{fig_pushing_the_curves_out_of_the_singularity}
 \caption{The subcurves of the $\gamma_i$ that are contained in $B(\sigma, \epsilon)$ are replaced by curves in $B(\sigma, \epsilon) \setminus \{\sigma\}$.}
 \end{center}
\end{figure}

 Now we do the same construction for the rest of the curves successively: For the construction of $\gamma_i'$ we seek a curve from $x_i^+$ to the prospective intersection point $z$ that does not leave $B(\sigma,\epsilon)$ and does not intersect the curves $\gamma_1',\ldots,\gamma_{i-1}'$. We find such a curve by taking any curve from $x_i^+$ to~$z$ in $B(\sigma,\epsilon) \setminus \{\sigma\}$ and instead of possibly crossing some~$\gamma_j'$ we follow~$\gamma_j'$ in a sufficiently small tubular neighbourhood without intersecting it until we reach~$z$.
 Then we define $\gamma_i'$ as the closed curve that consists of a curve from $z$ to $x_i^+$ as described, the subcurve of $\gamma_i$ from $x_i^+$ to $x_i^-$, and similarly a curve from $x_i^-$ to~$z$.
 Again, the curve $\delta_i$ is still a left-to-right curve of~$\gamma_i'$ with respect to $\gamma_1', \ldots, \gamma_{i-1}', \gamma_{i+1}, \ldots, \gamma_n$. 
 
 So we have a set of simple closed curves $\gamma_1',\ldots,\gamma_n'$ that are intersecting exactly in~$z$.
 Also, the set has left-to-right curves $\delta_1,\ldots, \delta_n$ so it is nonseparating by \autoref{lem_criterion_simple_closed_curves_nonseparating}.
\end{proof}
\end{lem}

\noindent
By the criterion in \autoref{lem_criterion_simple_closed_curves_nonseparating}, we can show that a set of curves is nonseparating but so far we do not have candidates of nonseparating curves for which we could use the criterion.
Therefore, we introduce a generalization of the well-known fact that saddle connections of finite translation surfaces are nonseparating.

For this, recall that the geodesic flow $F_\theta$ on a translation surface $(X, \mathcal{A})$ in a given direction $\theta \in S^1$ is an action of $\mathbb{R}$, defined on a subset of $X$. For a point $x\in X$ and a time $t\in \mathbb{R}$, the point $F_\theta(t, x)$ is the unique point in $X$ (if it exists) such that $x$ and $F_\theta (t, x)$ are connected by a geodesic segment of direction $\theta$ and length $t$. A geodesic flow $F_\theta$ is called \emph{recurrent} if for almost every point $x \in X$ it holds that, for every neighborhood $U \subseteq X$ of $x$ and every $t_0\in \mathbb{R}$, there exists a $t>t_0$ with $F_\theta(t,x) \in U$.

\begin{lem}[Saddle connections are nonseparating] \label{lem_saddle_connections_nonseparating}
 Let $(X,\mathcal{A})$ be a translation surface such that for two directions $\theta_1$, $\theta_2 \in S^1$ the geodesic flows $F_{\theta_1}$ and $F_{\theta_2}$ are recurrent.
 Furthermore, let $\gamma$ be a saddle connection starting and ending at the same singularity. Then $\gamma$ is nonseparating.
 
\begin{proof}
 Consider a geodesic segment $s\subseteq \im(\gamma)$ and the geodesic flow $F_\theta$ in a direction $\theta \in S^1$ that is transversal to the direction of $\gamma$ and such that the flow is recurrent.
 So there exists a point $x\in s$ that returns to $s$ under the flow $F_\theta$ after time $t_1$. In particular, there exists a time $t_0$ with $0<t_0 \leq t_1$ such that $F_\theta(x,t_0) \in \im(\gamma)$ for the first time. Additionally, as $X$ is orientable, $\gamma$ is geodesic, and $F_\theta$ is a geodesic flow, the curve $\delta \colon [0,t_0] \to X, \ t \mapsto F_\theta(x,t)$ is approaching $\im(\gamma)$ from the other side than it is leaving.
 Then for every $\epsilon > 0$ the curve $\delta_\epsilon \colon [\epsilon,t_0-\epsilon] \to X, \ t \mapsto F_\theta(x,t)$ is a curve in $X \setminus \im(\gamma)$. The curve $\delta_\epsilon$ or its reversed curve connects a point on the left side of $\gamma$ to a point on the right side of $\gamma$, so it is a left-to-right curve of $\gamma$. Hence, $\gamma$ is nonseparating.
\end{proof}
\end{lem}

\noindent
Note that there exist translation surfaces of infinite but also of finite area such that the geodesic flow is not recurrent for all but at most one direction. We describe such a translation surface of finite area in \autoref{exa_nested_cylinders} and note that it has in fact separating saddle connections.

We now have all ingredients to prove \autoref{prop_criterion_saddle_connections_infinite_genus}, which was stated at the beginning of the section.

\begin{proof}[Proof of \autoref{prop_criterion_saddle_connections_infinite_genus}]
 As saddle connections are curves that fulfill the conditions of \autoref{lem_nonseparating_curves_in_metric_completion}, we have a nonseparating set of $n$ simple closed curves in $X$ for every $n \geq 1$. By \autoref{def_genus}, this means that the genus of $X$ exceeds every number $n \in \mathbb{N}$, so $X$ has infinite genus.
\end{proof}

\section{Intersection of saddle connections} \label{sec_saddle_connections_intersections}

In the previous section, we showed that saddle connections are a good tool to prove infinite genus. However, these saddle connections cannot intersect arbitrarily and they need to have left-to-right curves.
Therefore, we introduce the immersion radius of saddle connections, which is related to the minimal length of intersecting saddle connections.

For translation surfaces with finitely many singularities which are all cone angle or infinite angle singularities, the lengths of the saddle connections are bounded from below. This is because for every singularity $\sigma$ there exists an $\epsilon_\sigma > 0$ such that $B(\sigma, \epsilon_\sigma) \setminus \{\sigma\}$ is locally flat. Therefore, the lengths of the saddle connections are bounded from below by $\min\{ \epsilon_\sigma : \sigma \text{ a singularity},\ B(\sigma, \epsilon_\sigma) \setminus \{\sigma\} \text{ locally flat}\}$.

For wild singularities, this argument does not work. In fact, for a wild singularity, the opposite is true.

\begin{prop}[Existence of short saddle connections]
 \label{prop_existence_short_saddle_connections}
 Let $(X,\mathcal{A})$ be a translation surface and $\sigma$ a wild singularity.
 Then for every $\epsilon > 0$, there exists a saddle connection connecting $\sigma$ to itself of length less than~$\epsilon$.

\begin{proof}
 As the singularities are discrete by \autoref{conv_discrete_singularities}, there exists an $\epsilon' > 0$ such that~$\sigma$ is the only singularity in $B(\sigma,\epsilon') \subseteq \overline{X}$.
 We distinguish between the following two cases.
 
 Case 1: There exists a rotational component of $\sigma$ which is bounded in at least one direction. This means that $\sigma$ itself is an obstacle to extending the rotational component in that direction. Hence, there exist curves from $\sigma$ to itself shorter than any given length, for instance shorter than~$\epsilon$. We now pass to the universal cover $\widetilde{X}$ of $X$ and consider a lift of such a short curve in the metric completion of $\widetilde{X}$. The metric completion of~$\widetilde{X}$ is a non-positively curved metric space and so, by \cite[Part II, Theorem 4.13]{bridson_haefliger_99}, there exists a geodesic in the metric completion of $\widetilde{X}$ which is homotopic to the lift of the original curve. As all additional points in the metric completion of $\widetilde{X}$ are projected to $\sigma$ in $\overline{X}$, the image of the chosen geodesic in $\overline{X}$ is either a saddle connection or it passes through $\sigma$ at least once. In the latter case, it contains a regular point and so a subsegment is in fact a saddle connection which then has length less than~$\epsilon$.

 Case 2: Every rotational component is unbounded.
 Assume that there exists a minimal length $\epsilon < \epsilon'$ of saddle connections from $\sigma$ to itself. Then for every linear approach $[\gamma]$ to $\sigma$, there exists a representative $\gamma \in \mathcal{L}^{\epsilon}(\sigma)$.
 Therefore, there exists a cyclic translation covering from $B(\sigma,\epsilon) \setminus \{\sigma\}$ to the once-punctured disk $B(0,\epsilon) \setminus \{0\} \subseteq \mathbb{R}^2$.
 This means that $\sigma$ is not a wild singularity.
\end{proof}
\end{prop}

\noindent
Even if there exist arbitrarily short saddle connections, a closed geodesic in $X$ cannot be intersected by arbitrarily short saddle connections as we will see in the discussion below. We can give a lower bound on the length of intersecting saddle connections by the \emph{immersion radius} of a closed geodesic.

The immersion radius of a point or of a curve is defined similarly to the well-known \emph{injectivity radius}. In contrast to the injectivity radius, we allow that the image of the disk that we map into $X$ overlaps itself.

\begin{definition}[Immersion radius] \label{def_immersion_radius}
 Let $(X,\mathcal{A})$ be a translation surface.
 
 \begin{enumerate}[series=immersionradiusitems]
  \item For a regular point $x\in X$, we define the \emph{immersion radius} $\ir(x)$ by
  \begin{equation*}
   \ir(x) \coloneqq d(x, \overline{X}\setminus X) \in (0,\infty] .
  \end{equation*}
  This is well-defined as $\overline{X}$ is a metric space and the set of singularities $\overline{X} \setminus X$ is discrete and hence closed in $\overline{X}$.
  Note that the open $\ir(x)$--neighbourhood of $x$ does not contain a singularity but its closure does (if $(X, \mathcal{A})$ has at least one singularity).
  
  \item For a curve $\gamma \colon [0,l] \to X$, we define the \emph{immersion radius} $\ir(\gamma)$ by
  \begin{equation*}
   \ir(\gamma) \coloneqq \inf\{ \ir( \gamma(t) ) : t \in [0,l] \} .
  \end{equation*}
  As the image of the curve $\gamma$ is compact, we can cover it by finitely many disks $B(\gamma(t_i), \ir(\gamma(t_i)))$ around points $\gamma(t_i)$. Then the set $\partial \big( \bigcup B(\gamma(t_i), \ir(\gamma(t_i))) \big)$ is compact and does not intersect $\im(\gamma)$, thus its distance to $\im(\gamma)$ is positive.
  As this distance is a lower bound for the immersion radius of the curve, this means that $\ir(\gamma) > 0$ still holds and again the open $\ir(\gamma)$--neighbourhood of $\im(\gamma)$ is locally flat.
  
  Note that every saddle connection that intersects~$\gamma$ has at least a length of $2 \ir(\gamma) > \ir(\gamma)$.
 \end{enumerate}
\end{definition}

\noindent
When generalizing the notion of immersion radius to saddle connections, we have to be careful as saddle connections are not curves in $X$ but in $\overline{X}$. Therefore we have to slightly modify the notion of immersion radius which compels us to restrict to so-called \emph{well-immersed} saddle connections. These are defined by considering the corresponding linear approaches.

\begin{rem}[Saddle connections and linear approaches] \label{rem_saddle_connections_linear_approaches}
 Recall that every saddle connection $\gamma \colon [0,l] \to \overline{X}$ is geodesic and has by definition an orientation. So the curve $\gamma_{|(0,l)}$ and its reversed curve define two linear approaches, one belonging to emanating from a singularity and one belonging to going into a singularity.
 We call the first linear approach $[\gamma_+]$ and the second one $[\gamma_-]$.
\end{rem}

\begin{definition}[Well-immersed saddle connections] \label{def_wellimmersed_saddle_connections}
 Let $(X,\mathcal{A})$ be a translation surface, $\sigma_+, \sigma_-$ two (not necessarily different) singularities and ${ \gamma \colon [0,l] \to \overline{X} }$ a saddle connection from $\sigma_+$ to $\sigma_-$. Let $[\gamma_+]$ and $[\gamma_-]$ be linear approaches as in \autoref{rem_saddle_connections_linear_approaches}.

 We say that $\gamma$ is a \emph{well-immersed} saddle connection if there exists an $\epsilon_+$ and an angular sector $( (0,\pi), 2\epsilon_+, i_{2\epsilon_+})$ with basepoint $\sigma_+$ such that $[\gamma_+] \in \im( f_{( (0,\pi), 2\epsilon_+, i_{2\epsilon_+})} )$ and the same is true for $[\gamma_-]$ with some $\epsilon_- > 0$ and angular sector $( (0,\pi), 2\epsilon_-, i_{2\epsilon_-})$ with basepoint~$\sigma_-$.
\end{definition}

\noindent
Note that for a rotational component of length strictly greater than $\pi$, there always exists an angular sector $( (0,\pi), \epsilon, i_{\epsilon})$ for some $\epsilon > 0$. This is true because the length of the longest representative of a linear approach varies lower-semicontinuously in a rotational component and hence has a minimum on a compact set of linear approaches (see \cite[Corollary 2.7 and the subsequent remark]{bowman_valdez_13}).
Hence, if for a saddle connection $\gamma$ the linear approaches $[\gamma_+]$ and~$[\gamma_-]$ are both contained in rotational components of length strictly greater than $\pi$ then $\gamma$ is well-immersed.

For well-immersed saddle connections, we can now also define the immersion radius.

\addtocounter{definition}{-3}
\begin{definition}[Immersion radius continued]\label{def_immersion_radius_saddle_connection}
\leavevmode\vspace{-\baselineskip} 
 \begin{enumerate}[resume=immersionradiusitems]
  \item Now let $\gamma \colon [0,l] \to \overline{X}$ be a well-immersed saddle connection from $\sigma_+$ to $\sigma_-$ and let $\gamma_+$ be the representative of $[\gamma_+]$ in $\mathcal{L}^{\epsilon_+}(\sigma_+)$.
  Because the length of $\overline{ [\gamma_+] }$ is at least $\pi$, the image of $\gamma_+$ is contained in a locally flat, open half-disk (see \autoref{fig_definition_immersion_radius_saddle_connection}).
  Suppose there exists a saddle connection that intersects $\gamma_+$. Then it has to start outside of the image of $i_{2\epsilon_+}$ and go through the half-annulus-like set
  \begin{equation*}
   \left( i_{2\epsilon_+} \circ f \right) \left( \left\{ z \in \mathbb{C} : \log \epsilon_+ \leq \Re(z) < \log 2 \epsilon_+, \Im(z) \in (0,\pi) \right\} \right)
  \end{equation*}
  to intersect the image of $\gamma_+$.
  This means that the intersecting saddle connection has at least length $\epsilon_+$.

  \begin{figure}
  \begin{center}
   \begin{tikzpicture}
    \draw[gray] (0,0) -- (-70:4cm);
    \draw[gray] (0,0) -- (120:4cm);
    \draw[gray] (-70:4cm) arc (-70:120:4cm);
    \draw[gray] (-20:4cm) node[right]{$2\epsilon_+$};
    \draw[gray] (-70:2cm) arc (-70:120:2cm);
    \draw[gray] (-40:2cm) node[right]{$\epsilon_+$};

    \draw[fill] (0,0) circle (1.3pt) node[left]{$\sigma_+$};
    \draw[very thick] (0,0) -- node[pos=0.7,below,inner sep=5pt]{$\gamma_+$} (35:2cm);
    \draw (35:2cm) -- (35:4.5cm);
  
    \draw[fill] (1.2,4.5) circle (1.3pt);
    \draw (1.2,4.5) -- (0.2,-1.6);
    \draw[fill] (0.2,-1.6) circle (1.3pt);
   \end{tikzpicture}  
  \label{fig_definition_immersion_radius_saddle_connection}
  \caption{Any saddle connection that intersects the well-immersed saddle connection~$\gamma$ close to its startpoint $\sigma_+$ has length at least $\epsilon_+$.}
  \end{center}
  \end{figure}
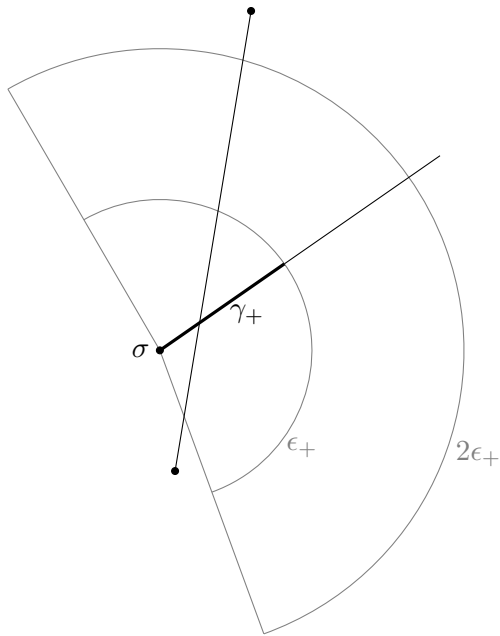

  A similar statement holds true for $\gamma_-$ and $\epsilon_-$.
  Furthermore, $\gamma_c \coloneqq \gamma_{| [\epsilon_+, l - \epsilon_-] }$ is a curve as in (ii) with a compact image and a well-defined immersion radius $\epsilon_c \coloneqq \ir(\gamma_c)>0$.

  Then $\min\{\epsilon_+, \epsilon_c, \epsilon_-\}$ is a lower bound for the length of saddle connections that intersect $\gamma$.
  We define the \emph{(generalized) immersion radius} $\ir(\gamma)$ by 
  \begin{align*}
   \ir(\gamma) \coloneqq \sup\{ \min\{\epsilon_+, \epsilon_c, \epsilon_-\} : & \ \epsilon_+ \text{ and } \epsilon_- \text{ small enough such that angular} \\ & \, \text{sectors as described above exist} \} .
  \end{align*}
 \end{enumerate}
\end{definition}
\addtocounter{definition}{2}

\noindent
In the case of a regular point $x\in X$, the term ``immersion radius'' is reasonable as an open disk of radius $\ir(x)$ can be immersed and the image under the immersion is a locally flat neighbourhood of the point~$x$. Similarly, for a closed geodesic $\gamma$, an open cylinder of height $2\cdot \ir(\gamma)$ and circumference the length of $\gamma$ can be immersed and the image is a locally flat tubular neighbourhood of the curve $\gamma$.

\begin{figure}
 \begin{center}
  \begin{tikzpicture}
   \draw[fill] (0,0) circle (1.2pt) node[left]{$\sigma$};
   \draw (0,0) -- node[below]{$\gamma$} (35:6cm);
   \draw[fill] (35:6cm) circle (1.2pt) node[right]{$\sigma$};

   \draw (0,0) -- (100:1.5cm);
   \draw (0,0) -- (-80:1.5cm);
   \draw[dotted] (0,0) -- (100:2cm);
   \draw[dotted] (0,0) -- (-100:2cm);
   \draw[dotted] (-100:2cm) arc (-100:100:2cm);

   \draw (35:6cm) -- +(140:1.4cm) -- (100:1.5cm);
   \draw (35:6cm) -- +(-40:1.4cm) -- (-80:1.5cm);
   \draw[dotted] (35:6cm) -- +(60:1.4cm);
   \draw[dotted] (35:6cm) -- +(-40:1.4cm) arc (-40:-300:1.4cm);
  \end{tikzpicture}
  \label{fig_trapezoidal_neighbourhood_of_saddle_connection}
  \caption{An example of a trapezoidal neighbourhood of a saddle connection $\gamma$.}
 \end{center}
\end{figure}
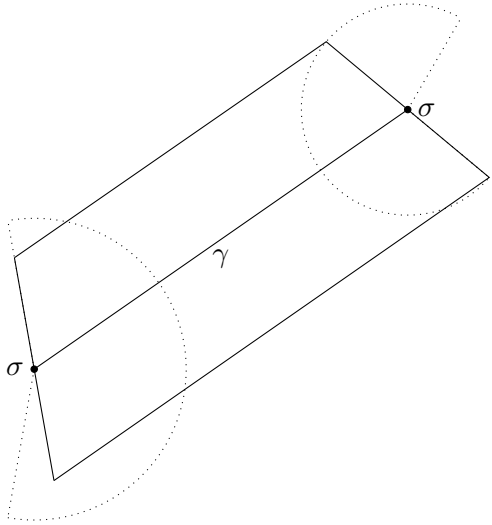

In the case of well-immersed saddle connections, we can immerse an open trapezoid of height $2 \cdot \ir(\gamma)$ such that the median has the same length as the saddle connection (see \autoref{fig_trapezoidal_neighbourhood_of_saddle_connection}).
The image forms a neighbourhood of the interior of the saddle connection and around the singularities we have the images of angular sectors of length $\pi$.

We finish this section by stating that the immersion radius varies continuously in the surface.

\begin{lem}[Immersion radius is continuous] \label{lem_immersion_radius_continuous}
 For a translation surface $(X, \mathcal{A})$, the map $\ir \colon X \to (0,\infty], \ x \mapsto \ir(x)$ is continuous.
 
\begin{proof}
 For all $x_1, x_2\in X$ and $\sigma \in \overline{X} \setminus X$, we have
  \begin{equation*}
  \ir(x_1) = d(x_1, \overline{X} \setminus X) \leq d(x_1,\sigma) \leq d(x_1,x_2) + d(x_2, \sigma) .
 \end{equation*}
 Therefore, it follows that
 \begin{equation*}
  \ir(x_1) \leq d(x_1, x_2) + \inf_{\sigma \in \overline{X} \setminus X} d(x_2, \sigma) = d(x_1, x_2) + \ir(x_2)
 \end{equation*}
 and interchanging $x_1$ and $x_2$ yields $| \ir(x_1) - \ir(x_2) | \leq d(x_1, x_2)$.
 So the map is Lipschitz continuous and hence continuous.
\end{proof}
\end{lem}

\section{Singularities that fulfill xossiness} \label{sec_xossiness}

The discussion on the immersion radius in the previous section shows that for well-immersed saddle connections, there exists a lower bound on the length of intersecting saddle connections.
To fulfill the conditions of \autoref{prop_criterion_saddle_connections_infinite_genus}, we are especially interested in translation surfaces for which there exist infinitely many saddle connections with such a lower bound on the length of intersecting saddle connections.
We say that the corresponding singularities fulfill \emph{xossiness} -- short for e\emph{x}istence \emph{o}f \emph{s}hort \emph{s}addle connections \emph{i}ntersected \emph{n}ot by \emph{e}ven \emph{s}horter \emph{s}addle connections.

\begin{definition}[Xossiness]
 Let $(X,\mathcal{A})$ be a translation surface and $\sigma$ a singularity.
 We say that $\sigma$ \emph{fulfills xossiness} if for every $\epsilon > 0$ there exists a saddle connection~$s$ that connects $\sigma$ to itself, that has length less than $\epsilon$, and is such that there exists a $\delta \coloneqq \delta (s) > 0$ such that no saddle connection of length less than $\delta$ intersects $s$.
\end{definition}

\noindent
Here and in the following, by saying ``two saddle connections do not intersect'' we mean that the images of the interiors of the saddle connections do not intersect. If we consider two saddle connections that connect the same singularity to itself then their images naturally have a common point, namely the singularity.

As the lengths of saddle connections starting in a fixed cone angle or infinite angle singularity are bounded away from $0$, cone angle and infinite angle singularities do not fulfill xossiness.
In the remainder of the section, we give two conditions for a wild singularity to fulfill xossiness.

\begin{prop}[Xossiness and the geodesic flow]\label{prop_xossiness_geodesic_flow}
 Let $(X,\mathcal{A})$ be a translation surface and $\sigma$ a wild singularity.
 Suppose that for a dense set of directions, for almost every point, the geodesic flow is defined for all time.
 Then $\sigma$ fulfills xossiness.

\begin{proof}
 Fix $\epsilon' > 0$ such that $\sigma$ is the only singularity in $B(\sigma,\epsilon') \subseteq \overline{X}$ and choose $\epsilon > 0$ with $\epsilon < \epsilon'$.
 By \autoref{prop_existence_short_saddle_connections}, there exists a saddle connection $s \colon [0,l] \to \overline{X}$ of length $l$ less than $\frac{\epsilon}{2}$.

 If there exists no saddle connection of length less than $\frac{\epsilon}{2}$ that intersects~$s$ then $s$ is a saddle connection as desired with $\delta(s) = \frac{\epsilon}{2}$.

 Now suppose there exists a saddle connection $s' \colon [0,l'] \to \overline{X}$ of length~$l'$ less than $\frac{\epsilon}{2}$ that intersects $s$ in $s(t)$ with $t\in (0,l)$.
 Choose a direction~$\theta$ such that for almost every point the geodesic flow $F_\theta$ is defined for all time and such that $\theta$ is so close to the direction of $s'$ that $F_\theta( s'(0), [0,l'])$ intersects~$s$ in $s(t_2)$ with $t_2 \in (0,t)$ and $F_\theta( s'(l'), [0,-l'])$ intersects $s$ in $s(t_3)$ with $t_3 \in (t,l)$ (see \autoref{fig_saddle_connection_intersected_by_shorter_saddle_connection}).
 Furthermore, choose $t_1\in (0,t_2)$ and $t_4 \in (t_3,l)$ such that $s(t_1)$ and $s(t_4)$ are points for which the geodesic flow $F_\theta$ is defined for all time, in particular backward and forward.
 
 Let $s([t_1,t_4])$ flow backward and forward under $F_\theta$ until it hits a singularity. This singularity is in both cases $\sigma$ because of the choice of $\epsilon < \epsilon'$.
 By this, we obtain an open parallelogram in $X$ that contains a singularity on two opposite edges but no singularities in the vertices as the vertices are images of $t_1$ or $t_4$ under $F_\theta$ (see \autoref{fig_saddle_connection_intersected_by_shorter_saddle_connection_swept_rectangle}).
 
 \begin{figure}
 \begin{center}
 \begin{tikzpicture}[scale=0.92]
  \draw[name path=s] (0,0) -- node[below, pos=0.85]{$s$} (6,2);
  \draw[fill] (0,0) circle (1pt) node[below]{$\sigma$};
  \draw[fill] (6,2) circle (1pt) node[below]{$\sigma$};

  \draw (2,-1) -- node[right, pos=0.85]{$s'$} (3,3);
  \draw[fill] (2,-1) circle (1pt);
  \draw[fill] (3,3) circle (1pt);

  \draw[color=gray, name path=forwardorbit] (2,-1) -- ++(60:4);
  \draw[color=gray, name path=backwardorbit] (3,3) -- ++(-120:4);

  \draw[fill, name intersections={of=s and backwardorbit, by={t_2}}] (t_2) circle (1pt) node[below]{$s(t_2)$};
  \draw[fill, name intersections={of=s and forwardorbit, by={t_3}}] (t_3) circle (1pt) node[below]{$s(t_3)$};
 \end{tikzpicture}
 \label{fig_saddle_connection_intersected_by_shorter_saddle_connection}
 \caption{A saddle connection $s$ from $\sigma$ to itself that is intersected by a saddle connection~$s'$ and two trajectories of the geodesic flow in a direction close to the direction of $s'$.}
 \end{center}
 \end{figure}
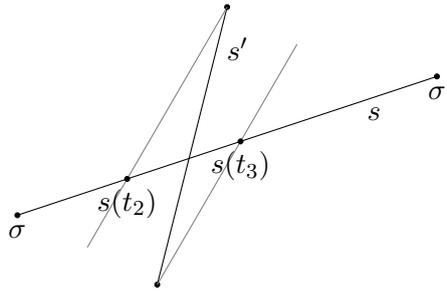

 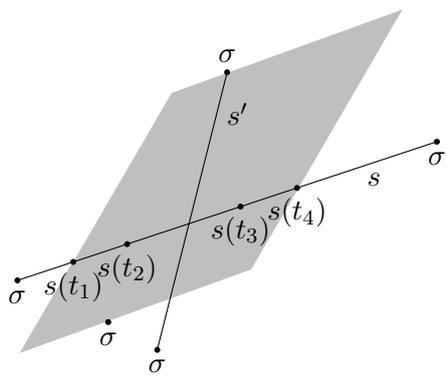
\begin{figure}
 \begin{center}
 \begin{tikzpicture}[scale=0.92]
  \draw[fill, color=gray!50] (0.8,0.267) -- ++(-120:1.5) -- ++(20:3.5) -- ++(60:4.3) -- ++(-160:3.5);

  \draw[name path=s] (0,0) -- node[below, pos=0.85]{$s$} (6,2);
  \draw[fill] (0,0) circle (1pt) node[below]{$\sigma$};
  \draw[fill] (6,2) circle (1pt) node[below]{$\sigma$};

  \draw[fill] (0.8,0.267) circle (1pt) node[below]{$s(t_1)$};
  \draw[fill] (4,1.333) circle (1pt) node[below]{$s(t_4)$};

  \draw (2,-1) -- node[right, pos=0.85]{$s'$} (3,3);
  \draw[fill] (2,-1) circle (1pt) node[below]{$\sigma$};
  \draw[fill] (3,3) circle (1pt) node[above]{$\sigma$};

  \path[color=gray, name path=forwardorbit] (2,-1) -- ++(60:4);
  \path[color=gray, name path=backwardorbit] (3,3) -- ++(-120:4);

  \draw[fill, name intersections={of=s and backwardorbit, by={t_2}}] (t_2) circle (1pt) node[below]{$s(t_2)$};
  \draw[fill, name intersections={of=s and forwardorbit, by={t_3}}] (t_3) circle (1pt) node[below]{$s(t_3)$};

  \draw[fill] (1.3,-0.6) circle (1pt) node[below]{$\sigma$};
 \end{tikzpicture}
 \label{fig_saddle_connection_intersected_by_shorter_saddle_connection_swept_rectangle}
 \caption{The segment $s([t_1,t_4)]$ is flowed forward and backward by $F_\theta$ until it hits a singularity for the first time.}
 \end{center}
\end{figure}

 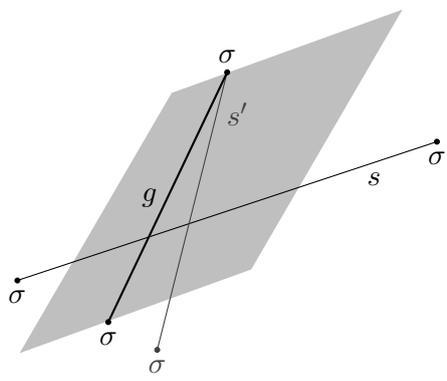
\begin{figure}
 \begin{center}
 \begin{tikzpicture}[scale=0.92]
  \draw[fill, color=gray!50] (0.8,0.267) -- ++(-120:1.5) -- ++(20:3.5) -- ++(60:4.3) -- ++(-160:3.5);

  \draw[name path=s] (0,0) -- node[below, pos=0.85]{$s$} (6,2);
  \draw[fill] (0,0) circle (1pt) node[below]{$\sigma$};
  \draw[fill] (6,2) circle (1pt) node[below]{$\sigma$};

  \draw[color=gray!50!black] (2,-1) -- node[right, pos=0.85]{$s'$} (3,3);
  \draw[fill, color=gray!50!black] (2,-1) circle (1pt) node[below]{$\sigma$};
  \draw[fill] (3,3) circle (1pt) node[above]{$\sigma$};

  \draw[fill] (1.3,-0.6) circle (1pt) node[below]{$\sigma$};

  \draw[thick] (1.3,-0.6) -- node[left]{$g$} (3,3);
 \end{tikzpicture}
 \label{fig_saddle_connection_intersected_by_shorter_saddle_connection_geodesic_in_rectangle}
 \caption{The new saddle connection $g$ cannot be intersected by saddle connections that are shorter than the distance of $g$ to the left edge and to the right edge of the parallelogram.}
 \end{center}
\end{figure}

 The geodesic $g$ in the parallelogram between the two appearances of $\sigma$ on the boundary as in \autoref{fig_saddle_connection_intersected_by_shorter_saddle_connection_geodesic_in_rectangle} is a saddle connection of length at most $l + l' \leq \epsilon$. The singularity $\sigma$ has a distinct distance to the left edge and to the right edge of the parallelogram and, as the interior of the parallelogram is locally flat, this distance is a lower bound for the immersion radius.
 Therefore, there exists a $\delta > 0$ such that no saddle connection of length less than $\delta$ can intersect $g$. Hence, $g$ is a saddle connection as desired.
\end{proof}
\end{prop}

\noindent
The second criterion uses a technical condition on the rotational components of a wild singularity.
The idea is again to find short well-immersed saddle connections.

\begin{prop}[Xossiness and rotational components] \label{prop_xossiness_wide_rotational_components}
 Let $(X,\mathcal{A})$ be a translation surface and $\sigma$ a wild singularity.
 Suppose that for every rotational component of $\sigma$ of length exactly $\pi$ there exists an angular sector $((0, \pi), \epsilon , i_\epsilon)$ such that the image of $f_{ ((0, \pi), \epsilon , i_\epsilon) }$ is contained in this rotational component.
 Then $\sigma$ fulfills xossiness.

\begin{proof}
 Fix $\epsilon' > 0$ such that $\sigma$ is the only singularity in $B(\sigma,\epsilon') \subseteq \overline{X}$ and choose $\epsilon > 0$ with $\epsilon < \epsilon'$.
 By \autoref{prop_existence_short_saddle_connections}, there exists a linear approach $[\gamma]$ for which the longest representative $\gamma$ has length $l < \frac{\epsilon}{2}$.

 For every $t\in (0,l)$, the immersion radius of $\gamma(t)$ is greater than $0$ but at most ${ d(\gamma(t),\sigma) \leq t }$. So we can define the \emph{immersion radius along $\gamma$} as the map
 \begin{equation*}
  \ir_{\gamma} \colon (0, l) \to (0, l), \ t \mapsto \ir(\gamma(t)) .
 \end{equation*}

 For every time $t\in (0,l)$, there exists a geodesic in $\overline{X}$ of length $\ir_\gamma(t)$ connecting $\gamma(t)$ to a singularity. Since $t + \ir_\gamma(t) < \epsilon'$, this singularity is again~$\sigma$.

 To prove that $\sigma$ fulfills xossiness, we show the existence of a time $t_0$ such that $\ir_\gamma(t_0)$ is realized by two different geodesics in $\overline{X}$. Then we can join the two occurrences of the singularity $\sigma$ at the endpoints of the geodesics in $B(\gamma(t_0),\ir_\gamma(t_0))$. The condition on the rotational components of length $\pi$ makes sure that we obtain a well-immersed saddle connection $s$ for which the immersion radius is defined as in \autoref{def_immersion_radius_saddle_connection} (iii). This means that there exists a lower bound on the length of saddle connections that intersect~$s$ in its interior.

 For every $t\in (0,l)$, there exists a locally flat disk $B(\gamma(t), \ir_\gamma(t))$.  We define the locally flat subset
 \begin{equation*}
  B \coloneqq \bigcup_{t\in (0,l)} B(\gamma(t), \ir_\gamma(t)) \subseteq X .
 \end{equation*}
 Furthermore, we consider a lift of $\gamma$ to the universal cover, together with lifts of $B(\gamma(t), \ir_\gamma(t) )$ for every $t\in (0,l)$. The union of the disks is a simply connected set and we develop it into the plane along the lift of $\gamma$ in the following way:
 Let $\widetilde{\gamma} \in \mathbb{R}^2$ be the holonomy vector of~$\gamma$ and define
 \begin{equation*}
  \widetilde{B} \coloneqq \bigcup_{t\in (0,l)} B(\widetilde{\gamma}(t), \ir_\gamma(t)) \subseteq \mathbb{R}^2 .
 \end{equation*}
 The set $\widetilde{B}$ is open and simply connected (see \autoref{fig_xossiness_developed_b}). Furthermore, there exists a map from $\widetilde{B}$ to $X$ which is a translation covering to $B$ and which can be extended to a map of the closure of $\widetilde{B}$ to $\overline{X}$.
 The preimages of $\sigma \in \overline{X}$ under this map are called \emph{representatives} of~$\sigma$. 
 Note that the singularity $\sigma$ is always the same on the boundary of $B \subseteq X$ while the representatives are not identified on the boundary of $\widetilde{B} \subseteq \mathbb{R}^2$.

 \begin{figure}
 \begin{center}
 \begin{tikzpicture}[rotate=20]
  \draw[thick] (0.3,0) circle (0.3cm);
  \draw[thick] (1,0) circle (0.5cm);
  \draw[thick] (2.5,0) circle (1.5cm);
  \draw[thick] (3.2,0) circle (1.1cm);
  \draw[thick] (4.1,0) circle (0.8cm);
  \draw[thick] (5.7,0) circle (1.7cm);
  \draw[thick] (6.5,0) circle (1.5cm);
  \fill[color=white] (0.3,0) circle (0.3cm);
  \fill[color=white] (1,0) circle (0.5cm);
  \fill[color=white] (2.5,0) circle (1.5cm);
  \fill[color=white] (3.2,0) circle (1.1cm);
  \fill[color=white] (4.1,0) circle (0.8cm);
  \fill[color=white] (5.7,0) circle (1.7cm);
  \fill[color=white] (6.5,0) circle (1.5cm);

  \fill (0,0) node[left] {$r$} circle (1pt);
  \fill (8,0) node[right] {$r'$} circle (1pt);
  \draw (0,0) -- node[above, pos=0.7] {$\widetilde{\gamma}$} (8,0);
 \end{tikzpicture}
 \label{fig_xossiness_developed_b}
 \caption{The open and simply connected, developed set $\widetilde{B} \subseteq \mathbb{R}^2$.}
 \end{center}
\end{figure}
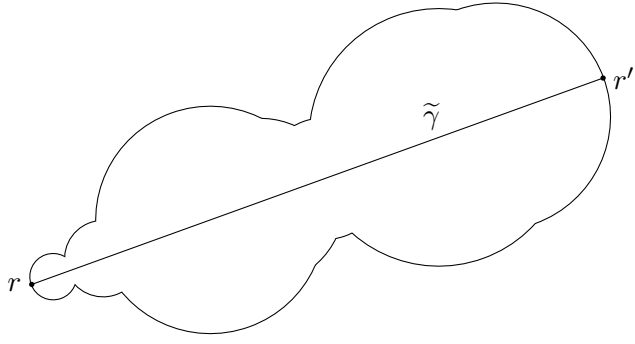
 
 Define the set $R \coloneqq \{r \in \partial \widetilde{B} : r \text{ is a representative of } \sigma\} \subseteq \mathbb{R}^2$.
 Every sequence in $R$ converging to a point $x\in \mathbb{R}^2$ corresponds to a sequence in~$\overline{X}$ where all elements are $\sigma$, so the limit is $\sigma$ and hence $x$ is a representative of~$\sigma$. Therefore, $R$ is a closed set in $\mathbb{R}^2$.

 For every representative $r \in R$, we define the set
 \begin{equation*}
  T_r \coloneqq \{ t \in (0, l) : d(\widetilde{\gamma}(t), r)=\ir_\gamma(t) \} .
 \end{equation*}
 Then the set $T_r$ is closed in $(0,l)$ and connected:

 \begin{itemize}
  \item Let $(t_n)_{n \in \mathbb{N}} \subseteq T_r$ be a sequence converging to a time $t \in (0, l)$. We have
  \begin{equation*}
   d(\widetilde{\gamma}(t_n), r) \leq d(\widetilde{\gamma}(t), r) + d(\widetilde{\gamma}(t), \widetilde{\gamma}(t_n)) = d(\widetilde{\gamma}(t), r) + |t-t_n|
  \end{equation*}
  and
  \begin{equation*}
   d(\widetilde{\gamma}(t_n), r) \geq d(\widetilde{\gamma}(t), r) - d(\widetilde{\gamma}(t), \widetilde{\gamma}(t_n)) = d(\widetilde{\gamma}(t), r) - |t-t_n| .
  \end{equation*}
  As $\ir \colon X \to (0,\infty]$ is continuous (see \autoref{lem_immersion_radius_continuous}), $\ir_\gamma$ is also continuous and we deduce
  \begin{equation*}
   \ir_\gamma(t) = \lim_{n\to \infty} \ir_\gamma(t_n) = \lim_{n \to \infty} d(\widetilde{\gamma}(t_n), r) = d(\widetilde{\gamma}(t), r) .
  \end{equation*}
  So $t$ is in $T_r$ and therefore $T_r$ is closed in $(0,l)$.

  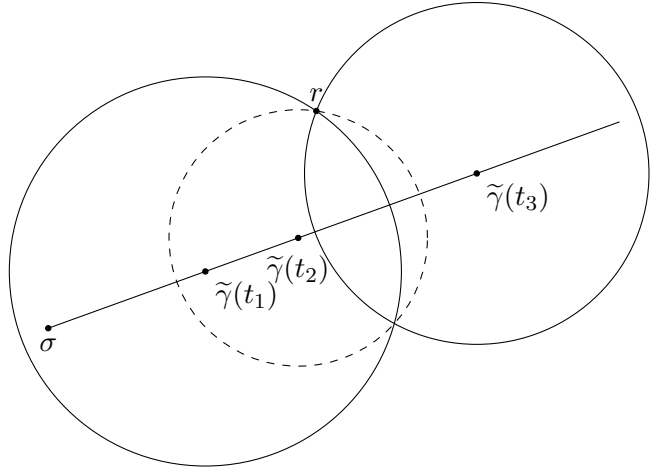
\begin{figure}
  \begin{center}
  \begin{tikzpicture}[rotate=20]
   \draw (0,0) -- (8,0);
   
   \draw[fill] (0,0) circle (1pt) node[below]{$\sigma$};
   \draw[fill] (2.2,0) circle (1pt) node[below right]{$\widetilde{\gamma}(t_1)$};
   \draw[fill] (3.5,0) circle (1pt) node[below=0.1cm]{$\widetilde{\gamma}(t_2)$};
   \draw[fill] (6,0) circle (1pt) node[below right]{$\widetilde{\gamma}(t_3)$};
   \draw[fill] (4.3,1.5) circle (1pt) node[above]{$r$};

   \node[draw] at (2.2,0) [circle through={(4.3,1.5)}] {};
   \node[draw] at (6,0) [circle through={(4.3,1.5)}] {};
   \node[draw, dashed] at (3.5,0) [circle through={(4.3,1.5)}] {};
  \end{tikzpicture}
  \label{fig_T_r_connected}
  \caption{The circle around $\widetilde{\gamma}(t_2)$ through $r$ is contained in the closure of the disks around $\widetilde{\gamma}(t_1)$ and $\widetilde{\gamma}(t_3)$ through $r$.}
  \end{center}
  \end{figure}

  \item For the proof of the connectedness consider $t_1$ and $t_3$ in $T_r$ and $t_2 \in (0, l)$ such that $t_1 < t_2 < t_3$.
  Then the circle around $\widetilde{\gamma}(t_2)$ through $r$ is contained in the closure of $B(\widetilde{\gamma}(t_1), \ir_\gamma(t_1)) \cup B(\widetilde{\gamma}(t_3), \ir_\gamma(t_3))$ (see \autoref{fig_T_r_connected}).
  This implies that for every $r' \in R$ with $d(\widetilde{\gamma}(t_2), r') < d(\widetilde{\gamma}(t_2), r)$ it holds as well that $d(\widetilde{\gamma}(t_1), r') < d(\widetilde{\gamma}(t_1), r)$ or $d(\widetilde{\gamma}(t_3), r') < d(\widetilde{\gamma}(t_3), r)$.
  Because we have chosen $t_1, t_3 \in T_r$, this is impossible and therefore we have $d(\widetilde{\gamma}(t_2), r) \leq d(\widetilde{\gamma}(t_2), r')$ for every $r' \in R$ and hence $t_2 \in T_r$.
 \end{itemize}

 We continue with a case-by-case analysis of $\partial \widetilde{B}$ and how it contains $R$:

 Case 1: There is an open, connected subset of $\partial \widetilde{B}$ which is disjoint from~$R$. Then there exists a closed connected subset $b$ of $\partial \widetilde{B}$ whose interior is disjoint from $R$ but whose endpoints (in a relative sense) are contained in $R$.
 We call these endpoints $r_1$ and $r_2$.

 To avoid technical subtleties, we now consider half-disks instead of disks $B(\gamma(t), \ir_{\gamma}(t))$ and slightly abuse notation. This means we only consider the connected components of $B(\widetilde{\gamma}(t), \ir_{\gamma}(t)) \setminus \im(\widetilde{\gamma})$ which are on the same side of $\widetilde{\gamma}$ as $r_1$ and $r_2$. Also, we only consider representatives in $R$ on the same side as $r_1$ and $r_2$, in $T_r$ we consider times $t$ where the geodesic from $\widetilde{\gamma}(t)$ to $r$ is the shortest on the chosen side, $\ir_\gamma(t)$ is the minimum of lengths of geodesics on the chosen side, and so on.

 In this sense we have that $T_{r_1} \cup T_{r_2}$ is connected:
 Choose $t_1 \in T_{r_1}$ and $t_2 \in T_{r_2}$ and assume $t_1 < t_2$.
 Furthermore, choose $t' \in (0,l)$ with $t_1 < t' < t_2$ (see \autoref{fig_xossiness_developed_b_connectedness}).
 Then we have for every $r \in R \setminus \{r_1,r_2\}$ on the same side of $\widetilde{\gamma}$ as $r_1$ and $r_2$ that the inequality
 $d(\widetilde{\gamma}(t'), r) \geq \min \{ d(\widetilde{\gamma}(t'), r_1), \, d(\widetilde{\gamma}(t'), r_2)\}$
 holds, by similar geometric arguments as in the proof of the connectedness of $T_r$.
 This means $t' \in T_{r_1} \cup T_{r_2}$ and $T_{r_1} \cup T_{r_2}$ is connected.

 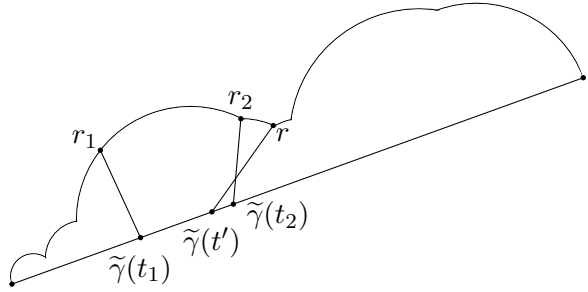
\begin{figure}[t]
 \begin{center}
 \begin{tikzpicture}[rotate=20]
  \begin{scope}
   \clip (0,0) rectangle (8,2);
   \draw[] (0.3,0) circle (0.3cm);
   \draw[] (1,0) circle (0.5cm);
   \draw[] (2.5,0) circle (1.5cm);
   \draw[] (3.2,0) circle (1.1cm);
   \draw[] (4.1,0) circle (0.8cm);
   \draw[] (5.7,0) circle (1.7cm);
   \draw[] (6.5,0) circle (1.5cm);
   \fill[color=white] (0.3,0) circle (0.29cm);
   \fill[color=white] (1,0) circle (0.49cm);
   \fill[color=white] (2.5,0) circle (1.49cm);
   \fill[color=white] (3.2,0) circle (1.09cm);
   \fill[color=white] (4.1,0) circle (0.79cm);
   \fill[color=white] (5.7,0) circle (1.69cm);
   \fill[color=white] (6.5,0) circle (1.49cm);
  \end{scope}

  \fill (1.7,1.27) node[above left=-0.1cm] {$r_1$} circle (1pt);
  \fill (3.58,1.03) node[above] {$r_2$} circle (1pt);
  \fill (3.95,0.8) node[below right=-0.1cm] {$r$} circle (1pt);
  \fill (1.8,0) node[below=0.1cm] {$\widetilde{\gamma}(t_1)$} circle (1pt);
  \draw (1.8,0) -- (1.7,1.27);
  \fill (3.1,0) circle (1pt);
  \draw (3.6,-0.35) node {$\widetilde{\gamma}(t_2)$};
  \draw (3.1,0) -- (3.58,1.03);
  \fill (2.8,0) node[below=0.07cm] {$\widetilde{\gamma}(t')$} circle (1pt);
  \draw (2.8,0) -- (3.95,0.8);
  
  \fill (0,0) circle (1pt);
  \fill (8,0) circle (1pt);
  \draw (0,0) -- (8,0);
 \end{tikzpicture}
 \label{fig_xossiness_developed_b_connectedness}
 \caption{If there exists no representative between $r_1$ and $r_2$ then $t'$ is contained in $T_{r_1} \cup T_{r_2}$.}
 \end{center}
\end{figure}
 
 As $T_{r_1}$ and $T_{r_2}$ are both closed and their union is connected, there exists a point $t_0 \in T_{r_1} \cap T_{r_2}$. In particular, the geodesics from $\widetilde{\gamma}(t_0)$ to $r_1$ and from $\widetilde{\gamma}(t_0)$ to $r_2$ in $\widetilde{B} \subseteq \mathbb{R}^2$ have the same length and the corresponding two geodesics from $\gamma(t_0)$ to $\sigma$ in $B\subseteq X$ are contained in $B(\gamma(t_0), \ir_\gamma(t_0))$ except for their endpoints.
 So we can join the endpoints of the two geodesics in $\overline{B} \subseteq \overline{X}$ and obtain a saddle connection from $\sigma$ to itself of length less than $2\cdot \ir_\gamma(t_0) \leq 2\cdot t_0 \leq 2l < \epsilon$ (see \autoref{fig_xossiness_developed_b_new_saddle_connection}).

 \begin{figure}[t]
 \begin{center}
 \begin{tikzpicture}[rotate=20]
  \begin{scope}
   \clip (0,0) rectangle (8,2);
   \draw[] (0.3,0) circle (0.3cm);
   \draw[] (1,0) circle (0.5cm);
   \draw[] (2.5,0) circle (1.5cm);
   \draw[] (3.2,0) circle (1.1cm);
   \draw[] (4.1,0) circle (0.8cm);
   \draw[] (5.7,0) circle (1.7cm);
   \draw[] (6.5,0) circle (1.5cm);
   \fill[color=white] (0.3,0) circle (0.29cm);
   \fill[color=white] (1,0) circle (0.49cm);
   \fill[color=white] (2.5,0) circle (1.49cm);
   \fill[color=white] (3.2,0) circle (1.09cm);
   \fill[color=white] (4.1,0) circle (0.79cm);
   \fill[color=white] (5.7,0) circle (1.69cm);
   \fill[color=white] (6.5,0) circle (1.49cm);
  \end{scope}

  \fill (1.7,1.27) node[above left=-0.1cm] {$r_1$} circle (1pt);
  \fill (3.58,1.03) node[above] {$r_2$} circle (1pt);
  \fill (2.5,0) circle (1pt);
  \draw (2.8,0) node[below=0.07cm] {$\widetilde{\gamma}(t_0)$};
  \draw (2.5,0) -- (1.7,1.27);
  \draw (2.5,0) -- (3.58,1.03);
  \draw[densely dashed] (1.7,1.27) -- (3.58,1.03);
  
  \fill (0,0) circle (1pt);
  \fill (8,0) circle (1pt);
  \draw (0,0) -- (8,0);
  
  \fill[color=gray!70!black] (1.8,0) node[below=0.05cm] {$\widetilde{\gamma}(t_1)$} circle (1pt);
  \fill[color=gray!70!black] (3.1,0) circle (1pt);
  \draw[color=gray!70!black] (3.6,-0.35) node {$\widetilde{\gamma}(t_2)$};
 \end{tikzpicture}
 \label{fig_xossiness_developed_b_new_saddle_connection}
 \caption{The dashed geodesic between $r_1$ and $r_2$ in $\widetilde{B}$ corresponds to a saddle connection in $B$.}
 \end{center}
\end{figure}
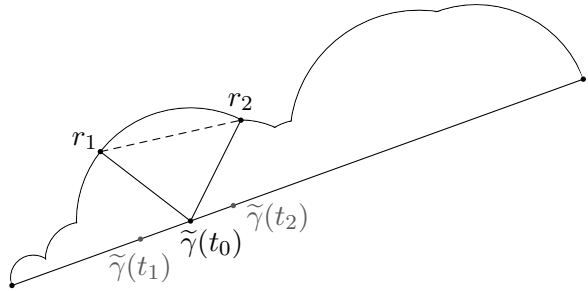
 
 The saddle connection from this construction is a chord of the half-disk $B(\gamma(t_0), \ir_\gamma(t_0))$, so the two corresponding linear approaches are contained in rotational components of length at least $\pi$.
 By the remark after \autoref{def_wellimmersed_saddle_connections} (for rotational components of length strictly greater than~$\pi$) and by assumption (for rotational components of length exactly~$\pi$) there exist angular sectors $((0, \pi), \epsilon_+, i_{\epsilon_+})$ and $((0, \pi), \epsilon_-, i_{\epsilon_-} )$ such that the linear approaches defined by the saddle connection are contained in one of the images of $f_{((0, \pi), \epsilon_+, i_{\epsilon_+})}$ and $f_{((0, \pi), \epsilon_-, i_{\epsilon_-} )}$.
 This means that the saddle connection is well-immersed and that there exists a lower bound on the length of saddle connections that intersect the resulting saddle connection in its interior.
 
 Case 2: Every open, connected subset of $\partial \widetilde{B}$ contains a representative~$r$. This means that $R$ is dense in $\partial \widetilde{B}$. As $R$ is a closed set, this implies $R = \partial \widetilde{B}$.
 When considering the situation in $X$, we have $\{\sigma\} = \partial B$, hence no point in $B$ can be connected to a point in~$X \setminus B$. As $X$ is connected, it follows $X=B$.
 Consider two parallel geodesics in $X$, which are close to each other and both end in $\sigma$, and a Cauchy sequence on each of these geodesics. Then the distance of the corresponding elements of the Cauchy sequences in $X$ is bounded away from $0$. This means that the limits of the Cauchy sequences are two different points but this is a contradiction as we only have one singularity in $\partial B = \{\sigma\}$.
 
 As Case 2 can never happen and, in Case 1, for each $\epsilon > 0$ we find a saddle connection as desired, the statement is proven.
\end{proof}
\end{prop}

\noindent
We have now seen in the proof how the peculiar condition on the existence of an angular sector $( (0, \pi), \epsilon, i_\epsilon)$ for every rotational component of length exactly $\pi$ is used. In particular, if such an angular sector does not exist then it is possible in Case 1 that short chords of the disk $B(\gamma(t_0), \ir_\gamma(t_0) )$ intersect the obtained saddle connection in its interior (cf.\ \autoref{fig_xossiness_developed_b_new_saddle_connection}).

\section{The main theorem, applications, and notes on generalizing}\label{sec_theorem_corollaries_generalizations}

In this section, we prove the main theorem that the existence of a wild singularity implies infinite genus under certain conditions. The key tool for the proof is \autoref{prop_criterion_saddle_connections_infinite_genus}. From the main theorem, we deduce three corollaries for classes of translation surfaces of particular interest.

Furthermore, we discuss whether the conditions in the main theorem are necessary and present an example of a translation surface that has a wild singularity and genus $0$.

\begin{thm}[Wild singularity implies infinite genus] \label{thm_wild_singularity_implies_infinite_genus}
 Let $(X,\mathcal{A})$ be a translation surface such that for two directions $\theta_1, \theta_2 \in S^1$ the geodesic flows $F_{\theta_1}$ and $F_{\theta_2}$ are recurrent.
 Furthermore, let $\sigma$ be a wild singularity of $(X,\mathcal{A})$ that fulfills xossiness.
 Then $X$ has infinite genus.
 
\begin{proof}
 According to the criterion in \autoref{prop_criterion_saddle_connections_infinite_genus}, to prove the statement we have to show that for every $n\geq 1$ that there exist $n$ saddle connections from~$\sigma$ to itself that intersect exactly in $\sigma$ and are such that the set has left-to-right curves. We do this by induction on $n$.
 
 For the base case $n=1$ we choose a saddle connection $\gamma_1$ such that there exists a lower bound $\epsilon_1 > 0$ on the length of intersecting saddle connections.
 Such a saddle connection exists by the assumption that $\sigma$ fulfills xossiness and, by \autoref{lem_saddle_connections_nonseparating}, it is nonseparating.

 For the inductive step assume that we have a nonseparating set of saddle connections $\{\gamma_1, \ldots, \gamma_{n-1} \}$ for which there exist lower bounds $\epsilon_1,\ldots, \epsilon_{n-1}$ on the lengths of intersecting saddle connections.
 In particular, for every $i\in \{1,\ldots, n-1\}$ there exists a left-to-right curve $\delta_i$, i.e.\ a curve in $X$ that connects the left side and the right side of $\gamma_i$ without intersecting any of the $\gamma_j$ for $j \in \{1,\ldots,n-1\}$.
 Let $\epsilon$ be the minimum of $\epsilon_1,\ldots,\epsilon_{n-1}$ and of the immersion radii of $\delta_1,\ldots,\delta_{n-1}$. As proven in \autoref{prop_existence_short_saddle_connections} there exists a saddle connection $\gamma_n$ from $\sigma$ to itself with length less than~$\epsilon$.
 Therefore, $\gamma_n$ does not intersect any of the curves $\gamma_1,\ldots,\gamma_{n-1}, \delta_1,\ldots,\delta_{n-1}$.
 Additionally, as $\sigma$ fulfills xossiness we can choose $\gamma_n$ such that there exists a lower bound $\epsilon_n > 0$ on the length of intersecting saddle connections.
 
 As $\gamma_n$ is nonseparating in $X$ (see again \autoref{lem_saddle_connections_nonseparating}) there exists a left-to-right curve $\delta_n'$ in $X\setminus \im(\gamma_n)$.
 If $\delta_n'$ does not intersect $\gamma_1,\ldots,\gamma_{n-1}$ then we can define $\delta_n \coloneqq \delta_n'$ and $\delta_n$ connects the left side of $\gamma_n$ to the right side of $\gamma_n$ in $X \setminus (\im(\gamma_1) \cup \ldots \cup \im(\gamma_{n}))$. Furthermore, none of the curves $\delta_1,\ldots,\delta_{n-1}$ intersects $\gamma_n$. Therefore we have that $X\setminus (\im(\gamma_1) \cup \ldots \cup \im(\gamma_n))$ is connected and the set of curves $\{ \gamma_1, \ldots, \gamma_n \}$ has left-to-right curves.

 If $\delta_n'$ intersects at least one of the curves $\gamma_1,\ldots,\gamma_{n-1}$ then we modify it in the following way. For every intersection with a curve $\gamma_i$ (without loss of generality, from the left of~$\gamma_i$) we choose a point $x_l$ on the left and a point $x_r$ on the right of $\gamma_i$ in $\im(\delta_n')$. Then we can replace the subcurve of $\delta_n'$ that intersects $\gamma_i$ by a curve in $N_l^\ast$ from $x_l$ to the startpoint of~$\delta_i$ concatenated with $\delta_i$ and concatenated with a curve in $N_r^\ast$ from the endpoint of $\delta_i$ to~$x_r$. By the induction hypothesis and by the choice of $\epsilon$, every $\delta_i$ for $i \in \{1,\ldots,n-1\}$ does not intersect any of the curves $\gamma_1,\ldots,\gamma_n$. Additionally, the new curve~$\delta_n$ is a left-to-right curve of~$\gamma_n$ with respect to $\gamma_1,\ldots,\gamma_{n-1}$.

 We have thus shown that for every $n\geq 1$, there exists a set of saddle connections of cardinality $n$ which has left-to-right curves.
 This implies by \autoref{prop_criterion_saddle_connections_infinite_genus} that $X$ has infinite genus.
\end{proof}
\end{thm}

\noindent
We now use the results from \autoref{prop_xossiness_geodesic_flow} and \autoref{prop_xossiness_wide_rotational_components} to give two possible conditions that make sure that \autoref{thm_wild_singularity_implies_infinite_genus} is applicable.

\begin{cor}[Wild singularity implies infinite genus] \label{cor_wild_singularity_implies_infinite_genus}
 Let $(X,\mathcal{A})$ be a translation surface and $\sigma$ a wild singularity. If one of the following two conditions is fulfilled then $X$ has infinite genus.
 \begin{enumerate}
  \item In a dense set of directions, the geodesic flow is recurrent.
  \item No rotational component of $\sigma$ has length $\pi$ and there exist two directions for which the geodesic flow is recurrent.
 \end{enumerate}
\end{cor}

\noindent
Note that the two properties in this corollary are not equivalent.
The second condition is more helpful for concrete given examples whereas the first condition can be used for more abstractly given classes of translation surfaces, such as parabolic translation surfaces and essentially finite translation surfaces.

\emph{Parabolic translation surfaces} are defined as having no Green's function. For them, it follows from Strebel's work \cite[Theorems 13.1 and 24.4]{strebel_84} that for all directions the geodesic flow is recurrent (cf.\ \cite[Remark 1]{trevino_12}).
Therefore we have the following corollary.

\begin{cor}[Wild parabolic translation surfaces have infinite genus]
 Let $(X, \mathcal{A})$ be a parabolic translation surface with a wild singularity. Then $X$ has infinite genus.
\end{cor}

\noindent
The other interesting class of translation surfaces was recently defined in \cite{rafi_randecker_18}.

\begin{definition}[Essentially finite translation surfaces]
 A translation surface $(X, \mathcal{A})$ is called \emph{essentially finite} if it fulfills the following three properties.
 \begin{enumerate}
  \item The surface $X$ has finite area.
  \item The set of singularities is discrete.
  \item Every singularity has countably many rotational components.
 \end{enumerate}
\end{definition}

\noindent
Examples are all finite translation surfaces as in \autoref{def_finite_translation_surfaces} and again the baker's map surface from \cite{chamanara_04} and the exponential surface from \cite{randecker_16}.

Note that every neighbourhood of an infinite angle singularity has infinite area so an essentially finite translation surface can only have cone angle singularities and wild singularities.
The number of these singularities is countable and, for every singularity, the number of linear approaches of a given direction is also countable. This implies that for a geodesic segment and a given transversal direction, the geodesic flow is defined for all but countably many points of the segment for all time. As the area is also finite, we can deduce from Poincaré recurrence that for all directions the geodesic flow on an essentially finite translation surface is recurrent.
Therefore, by \autoref{cor_wild_singularity_implies_infinite_genus}, an essentially finite translation surface with a wild singularity has infinite genus.

These observations help us to state a version of the Gauß--Bonnet formula for essentially finite translation surfaces.
The classical Gauß--Bonnet theorem relates the curvature of a surface $X$ to its Euler characteristic $\chi(X)$.
As the curvature of a translation surface is~$0$ everywhere, the formula for finite translation surfaces is given in a slightly different version in terms of the cone angles of the singularities. If $\ord(\sigma)$ is the multiplicity of a singularity $\sigma$ then it can be formulated as
\begin{equation*}
 - \chi(\overline{X}) = \sum\limits_{\sigma \in \overline{X}\setminus X} (\ord(\sigma_i) -1).
\end{equation*}
This is considered to be a formula in the spirit of the Gauß--Bonnet theorem as the singularities can roughly be seen as the points where the curvature of the closed surface~$\overline{X}$ is concentrated.

Together with the previous observation that essentially finite translation surfaces have no infinite angle singularities and with applying \autoref{cor_wild_singularity_implies_infinite_genus} to essentially finite translation surfaces, we can state the following version of the Gauß--Bonnet formula.

\begin{cor}[Gauß--Bonnet formula for essentially finite translation surfaces]
Let $(X, \mathcal{A})$ be an essentially finite translation surface. Then the genus $g$ of $X$ fulfills
\begin{equation*}
 2g-2 = \sum\limits_{\sigma \in \overline{X}\setminus X} (s_\sigma -1)
 \quad \text{where} \quad
 s_\sigma =
 \begin{cases}
  \ord(\sigma), & \sigma \text{ a cone angle singularity}, \\
  \infty, & \sigma \text{ a wild singularity}.
 \end{cases}
\end{equation*}
\end{cor}

\noindent
One of the key points of the proof of \autoref{thm_wild_singularity_implies_infinite_genus} is the assumption of recurrence of the geodesic flow.
By Poincaré recurrence one can deduce recurrence of the geodesic flow from the two weaker conditions that the flow is defined for a set of points of full measure for all time and that the area is finite.
We show that neither of the two conditions on its own works for the proof by considering the following two~examples.

\begin{exa}[Icicled surface] \label{exa_icicled_surface}
 Consider a half-open rectangle of height $2$ and width $1$. The left side is glued to the right side, the bottom and the top are excluded.

 For every $n \geq 1$, we consider a vertical segment starting at the bottom and a vertical segment starting at the top, at $\frac{i}{2^n}$ of length $\frac{1}{2^n}$ for every odd $i\in \{1,\ldots, 2^n-1\}$ (see \autoref{fig_icicled_surface}). We call the vertical segments~\emph{icicles}.

 \begin{figure}[t]
 \begin{center}
 \begin{tikzpicture}[scale=6]
  \draw[color=gray!30] (0,0) -- (1,0);
  \draw (1,0) -- (1,2);
  \draw[color=gray!30] (1,2) -- (0,2);
  \draw (0,2) -- (0,0);

  \foreach \zweihochn in {2,4,8,16,32,64} \foreach \i in {1,...,\zweihochn}
  {\draw (\i/\zweihochn,0) -- (\i/\zweihochn,1/\zweihochn);
  \draw (\i/\zweihochn,2) -- (\i/\zweihochn,2-1/\zweihochn);}
 \end{tikzpicture}
 \label{fig_icicled_surface}
 \caption{Vertical segments in the icicled surface.}
 \end{center}
 \end{figure}
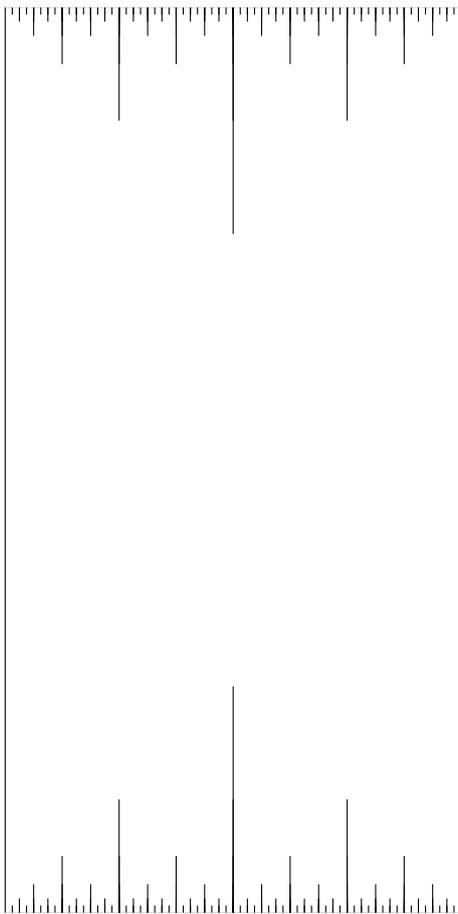

 Then we glue the segments as sketched in \autoref{fig_icicled_surface_gluings}. Note that no icicle on the top is glued to an icicle on the bottom. Formally, we can describe the gluings in the following way (starting with the icicles on the left part of the~top).
 \begin{itemize}
  \item For each side of the icicle at $\frac12$, we cut the segment again: first we cut it in half, then we cut the upper half into halves again, cut the upper quarter into halves again, \ldots \ So for every $n > 1$ we have a segment of length $\frac{1}{2^n}$ on each side of the icicle.
  \item The left side of the lower half of the icicle at $\frac12$ is glued to the right side of the icicle at $\frac14$.
  \item For every $n>2$ and every odd $i\in \{3,\ldots, 2^{n-1}-1\}$, the left side of the icicle at $\frac{i}{2^n}$ is glued to the right side of the icicle at $\frac{i-2}{2^n}$.
  \item For every $n>2$, the left side of the icicle at $\frac{1}{2^{n-1}}$ is cut into two segments of the same length. The lower part is glued to the right side of the icicle at $\frac{2^{n-1}-1}{2^n}$. The upper part is glued to the right side of the segment at $\frac12$ which has the correct length.
  \item We do the similar gluing for the right part of the top and also for the bottom.
 \end{itemize}
 
\begin{figure}
 \begin{center}
 \begin{tikzpicture}[scale=11, inner sep=2.5pt]
  \draw[dotted] (1,0.3) -- (1,0.4);
  \draw (1,0.4) -- (1,1);
  \draw[color=gray!30] (1,1) -- (0,1);
  \draw (0,1) -- (0,0.4);
  \draw[dotted] (0,0.4) -- (0,0.3);

  \foreach \zweihochn [evaluate=\zweihochn as \zweihochnminuseins using \zweihochn-1] in {2,4,8,16}
  \foreach \i in {1,...,\zweihochnminuseins}
  {\draw (\i/\zweihochn,1) -- (\i/\zweihochn,1-1/\zweihochn);}

  \fill (0.5,0.75) circle (0.1pt);
  \fill (0.5,0.875) circle (0.1pt);
  \fill (0.5,0.9375) circle (0.1pt);
  \fill (0.25,0.875) -- +(0,0.1pt) arc(90:270:0.1pt) -- (0.25,0.875);
  \fill (0.125,0.9375) -- +(0,0.1pt) arc(90:270:0.1pt) -- (0.125,0.9375);
  \fill (0.0625,0.96875) -- +(0,0.1pt) arc(90:270:0.1pt) -- (0.0625,0.96875);

  \path (0.5,0.9375) -- node[right] {$z_2$} (0.5,0.875) -- node[right] {$z_1$} (0.5,0.75) -- node[left] {$a$} (0.5,0.5);
  \path (0.25,1) -- node[right] {$a$} (0.25,0.75);
  \path (0.25,1) -- node[left, inner sep=1.5pt] {$z_1$} (0.25,0.875) -- node[left] {$b$} (0.25,0.75);
  \path (0.375,1) -- node[right] {$b$} node[left, text height=height("b")] {$c$} (0.375,0.875);
  \path (0.125,1) -- node[right] {$c$} (0.125,0.875);
  \path (0.125,1) -- node[left, inner sep=1.5pt] {$z_2$} (0.125,0.9375) -- node[left] {$d$} (0.125,0.875);
  \path (0.4375,1) -- node[right] {$d$} node[left,  text height=height("d")] {$e$} (0.4375,0.9375);
  \path (0.3125,1) -- node[right] {$e$} node[left] {$f$} (0.3125,0.9375);
  \path (0.1875,1) -- node[right] {$f$}(0.1875,0.9375);
 \end{tikzpicture}
  
 \label{fig_icicled_surface_gluings}
 \caption{Gluings for the icicled surface: segments with the same letters are glued.}
 \end{center}
\end{figure}
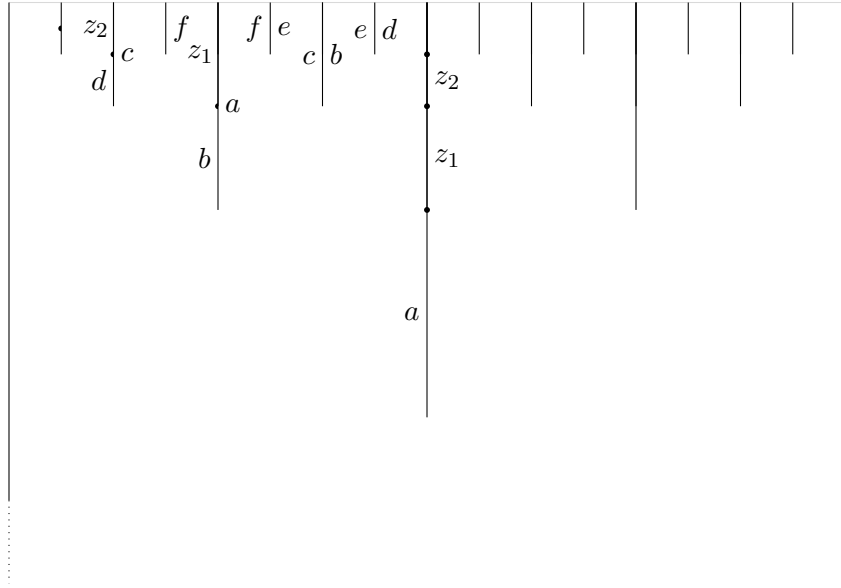

The resulting translation surface $(X, \mathcal{A})$ is called \emph{icicled surface} and has the following properties.
 \begin{enumerate}
  \item There exist exactly two singularities: all the tips of the icicles on the top are identified by the definition of the gluings. We call the corresponding singularity~$\sigma_{top}$. Now consider a nondyadic point $p$ in the top boundary of the rectangle, i.e.\ a point where no icicle starts. There exists a sequence of icicles such that the tips of the icicles converge to $p$, seen as points in $\mathbb{R}^2$ without gluings. Therefore, the distance from $\sigma_{top}$ to $p$ is $0$ and $p=\sigma_{top}$ in $\overline{X}$. The same argument works for the dyadic points on the boundary where an icicle starts. Note that the points where we cut some icicles into more segments are glued to points on the boundary, so these are also identified with $\sigma_{top}$.
  
  Also, the same reasoning holds for the tips of the icicles on the bottom and the points in the bottom boundary. So we have only two singularities $\sigma_{top}$ and $\sigma_{bottom}$. In particular, the set of singularities is discrete.

  Note that both of the singularities have one rotational component of infinite length and uncountably many rotational components of finite length, so the icicled surface is not essentially finite.

  \item Every icicle defines a saddle connection or a chain of saddle connections.
  However, for most of these saddle connections defined by icicles, there also exist arbitrarily short saddle connections close to the top or bottom, intersecting them.
  On the other hand, we can define saddle connections from the tip of the icicle at $\frac{1}{2^n}$ to the tip of the icicle at $\frac{1}{2^{n+1}}$ for every $n\geq 1$ (see \autoref{fig_icicled_surface_saddle_connections}). The length of the $n$th such saddle connection is $\smash{\frac{\sqrt{2}}{2^{n+1}}}$ and no saddle connection of length smaller than $\smash{\frac{\sqrt{2}}{2^{n+1}}}$ can intersect it. This means that both singularities fulfill xossiness.
  
  \begin{figure}
  \begin{center}
   \begin{tikzpicture}[scale=6]
    \draw[dotted] (1,0.3) -- (1,0.4);
    \draw (1,0.4) -- (1,1);
    \draw[color=gray!30] (1,1) -- (0,1);
    \draw (0,1) -- (0,0.4);
    \draw[dotted] (0,0.4) -- (0,0.3);

    \foreach \zweihochn [evaluate=\zweihochn as \zweihochnminuseins using \zweihochn-1] in {2,4,8,16,32,64}
    \foreach \i in {1,...,\zweihochnminuseins}
    {\draw (\i/\zweihochn,1) -- (\i/\zweihochn,1-1/\zweihochn);}

    \fill (0.5,0.75) circle (0.17pt);
    \fill (0.5,0.875) circle (0.17pt);
    \fill (0.5,0.9375) circle (0.17pt);
    \fill (0.25,0.875) -- +(0,0.17pt) arc(90:270:0.2pt) -- (0.25,0.875);
    \fill (0.125,0.9375) -- +(0,0.17pt) arc(90:270:0.2pt) -- (0.125,0.9375);
    \fill (0.0625,0.96875) -- +(0,0.17pt) arc(90:270:0.2pt) -- (0.0625,0.99875);

    \draw[very thick] (0.5,0.5) -- (0.25,0.75);
    \draw[very thick] (0.125,0.875) -- (0.0625,0.9375);
   \end{tikzpicture}
   \label{fig_icicled_surface_saddle_connections}
   \caption{The first and the third saddle connection of the set described in \autoref{exa_icicled_surface}~(ii) to show xossiness.}
   \end{center}
  \end{figure}
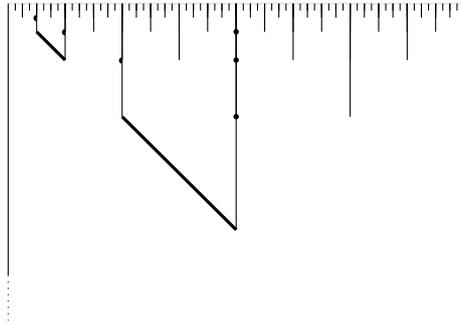
  
  \item For a regular point in $X$, the trajectory under a geodesic flow is only defined until it hits the tip of an icicle or the top or bottom of the rectangle. Hence, there exists no point in $X$ such that the geodesic flow in the vertical direction is defined for all time.
  
  Let $\theta$ be a direction in $(0,\pi)$. Suppose that the geodesic flow $F_\theta$ is defined for all time for a set of points of full measure. Consider a closed horizontal geodesic $g$ in the middle of the surface and a tubular neighbourhood $N$ of $g$ not intersecting any icicles. Then $F_\theta$ is defined on a subset of $N$ of full measure for all time.
  From Poincaré recurrence we can deduce that there exists a point $x \in N$ and a time $t_x>0$ such that $F_\theta(x, t_x) \in N$. This means that $F_\theta(x, t_x-\epsilon)$ is contained in the lower part of the surface for an $\epsilon > 0$. But this is impossible as there does not exist a way to reach the lower part of the surface from the upper part of the surface without intersecting $N$.

  Hence for all directions except for the horizontal one, there is no set of points of full measure for which the geodesic flow is defined for all time.

  \item Because all but the horizontal geodesic flow are not recurrent, we cannot use the arguments of \autoref{thm_wild_singularity_implies_infinite_genus} to show infinite genus.
  On the other hand, we can check by a sharp look that every icicle is defining at least one saddle connection which has a left-to-right curve. Moreover, every set of saddle connections defined by icicles has left-to-right curves and the number of icicles is not bounded. Therefore, for every $n\geq 1$ there exists a set of $n$ saddle connections that has left-to-right curves. With this we can show directly by \autoref{prop_criterion_saddle_connections_infinite_genus} instead of \autoref{thm_wild_singularity_implies_infinite_genus} that $X$ has infinite genus.
 \end{enumerate}
\end{exa}

\noindent
We want to emphasize that in this example, the top and bottom boundary is in some sense ``trapping'' the whole geodesic flow. This kind of dynamical behaviour was so far only known for translation surfaces with a continuum of singularities, for instance for the open disk (see the remark before \autoref{conv_discrete_singularities}).
However, we explicitly excluded this kind of example by requiring that our translation surfaces have discrete singularities.
There was the expectation that discreteness of the set of singularities should imply good dynamical properties.
However, when starting with a continuum of singularities such as the open rectangle in \autoref{exa_icicled_surface} and introducing additional gluings, we can maintain the bad dynamical properties and identify the continuum of singularities to one point. So the icicled surface is morally derived from a translation surface with non-discrete singularities but formally it shows that there exist translation surfaces with a discrete set of singularities so that for at most one direction $\theta$ the geodesic flow $F_\theta$ is defined for almost every point for all time. In particular, for such an example it is not possible to apply Poincaré recurrence to conclude from the finiteness of the area that a geodesic flow is recurrent.

The recurrence of the flow is needed for the proof of \autoref{thm_wild_singularity_implies_infinite_genus} as it is an assumption in \autoref{lem_saddle_connections_nonseparating}, which states the crucial ingredient of the proof that saddle connections are nonseparating under the condition on recurrence.
In fact, there are translation surfaces (with discrete singularities) that have separating saddle connections.
For example, the horizontal saddle connection of the icicled surface that connects the tip of the longest icicle to itself is separating. However, as was indicated in \autoref{exa_icicled_surface} (v) there exist other saddle connections that are nonseparating. This is not necessarily the case, as we will see in the next example, which was worked out together with Pat Hooper.

\begin{exa}[Nested cylinders] \label{exa_nested_cylinders}
 Consider a Euclidean half-plane with a distinguished midline. We cut vertical slits of infinite length in the half-plane from the midline upward, starting from $\sum_{i=1}^{n}  \frac{1}{i}$ for $n$ odd.
 Additionally, we cut vertical slits of infinite length from the midline downward, starting from $\sum_{i=1}^{n}  \frac{1}{i}$ for $n$ even.
 Now we glue the right side of a slit to the left side of the slit which is next to the right and the left side of the slit to the right side of the slit which is next on the left (see \autoref{fig_wild_singularity_genus_0}).
 By this construction, we obtain half-cylinders with smaller and smaller circumferences that are glued in a nested way.
 
 \begin{figure}
 \begin{center}
  \begin{tikzpicture}[scale=0.42]
   \draw (0,0) -- (0,5) -- node[right, text height=height("b")] {$a$} (0,9);
   \draw (6,0) .. controls +(110:1cm) and +(-90:2cm) .. (5.5,5) -- node[left, text height=height("b")] {$a$} (5.5,9);
   \draw (6,0) .. controls +(70:1cm) and +(-90:2cm) .. (6.5,5) -- node[right] {$b$} (6.5,9);
   \draw (13,0) .. controls +(110:0.8cm) and +(-90:2cm) .. (12.6,5) -- node[left] {$b$} (12.6,9);
   \draw (13,0) .. controls +(70:0.8cm) and +(-90:2cm) .. (13.4,5) -- node[right, text height=height("b")] {$c$} (13.4,9);
   \draw (17.4,0) .. controls +(110:0.6cm) and +(-90:2cm) .. (17.1,5) -- node[left, text height=height("b")] {$c$} (17.1,9);
   \draw (17.4,0) .. controls +(70:0.6cm) and +(-90:2cm) .. (17.7,5) -- node[right] {} (17.7,9);
   \draw (20.614,0) .. controls +(110:0.6cm) and +(-90:2cm) .. (20.414,5) -- node[left] {} (20.414,9);
   \draw (20.614,0) .. controls +(70:0.6cm) and +(-90:2cm) .. (20.814,5) -- node[right] {} (20.814,9);

   \draw (0,0) -- (0,-5) -- node[right] {$d$} (0,-9);
   \draw (10,0) .. controls +(250:0.8cm) and +(90:2cm) .. (9.6,-5) -- node[left] {$d$} (9.6,-9);
   \draw (10,0) .. controls +(290:0.8cm) and +(90:2cm) .. (10.4,-5) -- node[right, text height=height("b")] {$e$} (10.4,-9);
   \draw (15.4,0) .. controls +(250:0.6cm) and +(90:2cm) .. (15.1,-5) -- node[left, text height=height("b")] {$e$} (15.1,-9);
   \draw (15.4,0) .. controls +(290:0.6cm) and +(90:2cm) .. (15.7,-5) -- node[right] {$f$} (15.7,-9);
   \draw (19.114,0) .. controls +(250:0.6cm) and +(90:2cm) .. (18.914,-5) -- node[left] {$f$} (18.914,-9);
   \draw (19.114,0) .. controls +(290:0.6cm) and +(90:2cm) .. (19.314,-5) -- node[right] {} (19.314,-9);

   \fill (0,0) circle (2.5pt);

   \draw (22.5,5) node {$\cdots$};
   \draw (21,-5) node {$\cdots$};
  \end{tikzpicture}
 \label{fig_wild_singularity_genus_0}
 \caption{The nested cylinders example has a wild singularity and genus $0$.}
 \end{center}
 \end{figure}
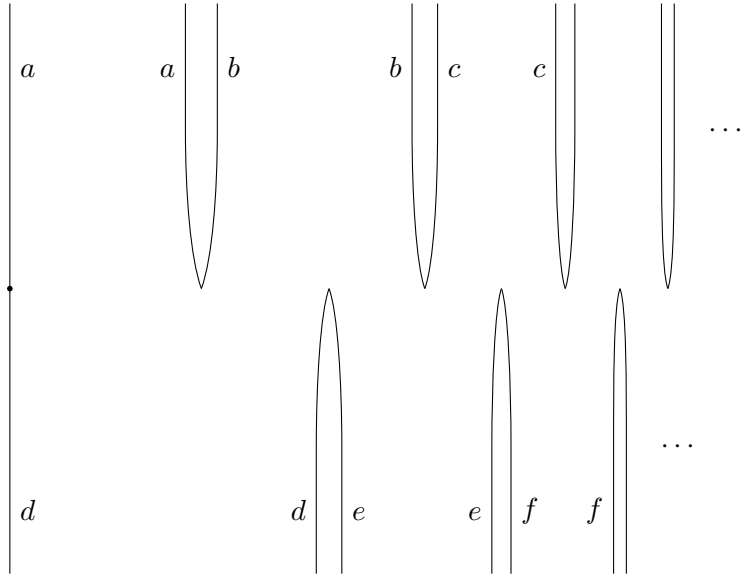
 
 The resulting translation surface has infinite area, genus $0$, and exactly one singularity. The singularity is wild as the distance between the startpoints of the slits is going to $0$, i.e.\ there exists no cyclic translation covering from a punctured neighbourhood of $\sigma$ to a once-punctured disk in $\mathbb{R}^2$. This singularity has exactly one rotational component which is isometric to $\mathbb{R}$.
 
 For every direction $\theta$, the geodesic flow $F_\theta$ is defined for almost every point for all time. However, recurrence only occurs in the horizontal direction. In this example, Poincaré recurrence is not applicable because the area is not finite.
 
 In particular, all saddle connections are horizontal and all of them are separating.
\end{exa}

\noindent
The last example indicates that the statement in \autoref{thm_wild_singularity_implies_infinite_genus} is wrong if we give up the condition on recurrence in two directions.
However, the example does not destroy the prospect of replacing the recurrence condition by a weaker condition like finite area.

\bibliographystyle{amsalpha}
\bibliography{/home/anja/Documents/Mathematik/literature/BibTex/Literatur}

\end{document}